\documentclass[leqno,10pt]{article}
\usepackage{amsmath}
\usepackage{amscd,amssymb}
\usepackage{verbatim}
\usepackage[all]{xy}

\newtheorem{Lemma}{Lemma}[section]
\newtheorem{Th}[Lemma]{Theorem}
\newtheorem{Prop}[Lemma]{Proposition}

\newtheorem{Cor}[Lemma]{Corollary}
\newtheorem{Def}[Lemma]{Definition}
\newtheorem{Ex}[Lemma]{Example}
\newtheorem{Fc}[Lemma]{Fact}
\newtheorem{Exs}[Lemma]{Examples}

\newtheorem{Remark}[Lemma]{Remark}
\newenvironment{Proof}{{\sc Proof.}\ }{~\rule{1ex}{1ex}\vspace{0.5truecm}}

\newcommand{\Hom}{\mbox{\rm Hom}}

\newcommand{\Ext}{\mbox{\rm Ext}}
\newcommand{\Tor}{\mbox{\rm Tor}}

\newcommand{\card}{\mbox{\rm card}}

\newcommand{\Add}{\mbox{\rm Add}}

\newcommand{\B}{\mathcal{B}}

\newcommand{\N}{\mathbb N}

\newcommand{\C}{\mathcal{C}}
\newcommand{\D}{\mathcal{D}}

\newcommand{\Q}{\mathcal{Q}}
\newcommand{\A}{\mathcal{A}}
\newcommand{\Y}{\mathcal{Y}}
\newcommand{\F}{\mathcal{F}}
\newcommand{\X}{\mathcal{X}}

\newcommand{\Mod}{\mbox{\rm Mod-}}
\newcommand{\lmod}{\mbox{\rm -mod}}
\newcommand{\lMod}{\mbox{\rm -Mod}}

\newcommand{\LMod}{\mbox{\rm -Mod}}
\newcommand{\rmod}{\mbox{\rm mod-}}
\newcommand{\mk}[1]{\mathrm{ mod}_{#1}\mbox{-}R}
\newcommand{\qk}[1]{\mathrm{ mod}_{#1}\mbox{-}Q}
\newcommand{\abs}[1]{\mid {#1}\vert}

\newcommand{\PP}{\ensuremath{\mathcal{P}}}
\newcommand{\ES}{\ensuremath{\mathcal{S}}}

\newcommand{\TFcal}{\ensuremath{\mathcal{TF}}}

\newcommand{\m}{\ensuremath{\mathbf{m}}}

\begin{document}
\title {Cotorsion pairs generated by modules of bounded projective dimension}

\author{Silvana Bazzoni\thanks{Supported by  grant SAB2005-0139 of the Secretar\'{i}a de Estado de Univesidades e
Investigaci\'{o}n del Ministerio de Educaci\'{o}n y Ciencia de
Espa\~na.
 Partially supported by  Universit\`{a} di Padova (Progetto di Ateneo
CDPA048343 ``Decomposition and tilting theory in modules, derived and
cluster categories'').} \\
Dipartimento di Matematica Pura e Applicata, \\ Universit\`a di
Padova\\ Via Trieste 63, 35121 Padova, Italy \\ e-mail:
bazzoni@math.unipd.it \and Dolors Herbera\protect\thanks {Partially
supported by MEC-DGESIC (Spain) through Project MTM2005-00934, and
by the Comissionat Per Universitats i Recerca de la Generalitat de
Catalunya through Project 2005SGR00206.  \protect\newline While this
paper was written, both authors were working within the Research
Programme
 on Discrete and Continuous methods of Ring Theory at the CRM, Barcelona
 (Spain). They thank their host for its hospitality.
  \protect\newline 2000 Mathematics
Subject Classification. Primary: 16D90; 16E30; Secondary:  16G99.}
\\ Departament de Matem\`atiques, \\
Universitat Aut\`onoma de Barcelona, \\ 08193 Bellaterra
(Barcelona), Spain\\ e-mail: dolors@mat.uab.cat }

\date{\phantom{ciao}}

\maketitle
\begin{abstract} We apply the theory of cotorsion pairs to  study closure properties of classes of
modules with finite projective dimension with respect to direct
limit operations and  to   filtrations.

We also prove that if the ring is an order in an $\aleph
_0$-noetherian ring $Q$ of small finitistic dimension $0$, then  the
cotorsion pair generated by the modules of projective dimension at
most one is of finite type if and only if $Q$ has big finitistic
dimension $0$. This applies, for example, to semiprime Goldie rings
and Cohen Macaulay noetherian commutative rings.\end{abstract}

\section{Introduction}

In this paper we apply the theory of cotorsion pairs to study
classes of modules with finite projective dimension. The first
insight in this direction was made in \cite{AT} (see also
\cite[Chapter~7]{T2}). Our approach takes advantage, and it is based
on, the recent developments in the area that had led, for example,
to show that all tilting modules are of finite type \cite{BET},
\cite{BH}, \cite{ST1}, \cite{BS} and to solve the Baer splitting
problem raised by Kaplansky in 1962 \cite{abh}.

For a ring $R$, let $\PP_n$ be the class of right $R$-modules of
projective dimension at most $n$. Denote  by $\rmod R$ the resolving
class of right $R$-modules having a projective resolution consisting
of finitely generated projective modules. We set $\rmod R\ \cap\ \PP
_n : = \PP_n(\rmod R)$.

A possible approach to understand the structure of the modules in
$\PP _n$ in terms of modules in $ \PP_n(\rmod R)$, is to determine
whether they belong to the direct limit closure of  $ \PP_n(\rmod
R)$. However, as direct limits do not commute with the
$\mathrm{Ext}$ functor, it is also convenient to   turn the
attention towards the smaller class of modules filtered by modules in $\PP_n(\rmod
R)$ or, even better, towards direct summands of such modules. (See Fact~\ref{F:Jan-description}.)


From the general theory of cotorsion pairs it follows that the
modules in $\PP_n$ are direct summands of $\PP_n(\rmod R)$-filtered
modules if and only if the cotorsion pair $(\PP_n, \PP_n^\perp)$ is
of finite type, that is,  if and only if $\PP_n^\perp= \PP_n(\rmod
R)^\perp$. ($^\perp$ denotes the $\Ext$-orthogonal, see \S
~\ref{preliminaries} for unexplained terms and notation).

We summarize these results, as well as the relation between
filtrations and direct limits in Proposition~\ref{P:finite}. We give
a new insight on this interaction in Theorem~\ref{T:lim-in P_n},
where we show that if $(\PP_n, \PP_n^\perp)$ is of finite type, then
any module in $\PP_{n+1}$ is a direct limit of modules in
$\PP_{n+1}(\rmod R)$. On the other hand we also prove that the
finite type of the cotorsion pair $(\PP_n, \PP_n^\perp)$ for some
$n\ge 1$ implies strong coherency/noetherianity conditions on the
class $\PP_n$, see Corollary~\ref{C:coherence}.

The basic idea to show the finite type of the cotorsion pair
$(\PP_n, \PP_n^\perp)$, is patterned in the method used to prove
that tilting classes are of finite type. This means to follow a
two-step procedure: First to show that $(\PP_n, \PP_n^\perp)$  is a
cotorsion pair of countable type, and then conclude the finite type
by proving that the $\Ext$-orthogonal of the countably presented
modules coincides with the $\Ext$-orthogonal of the finitely
presented ones.

After Raynaud a Gruson \cite{RG}, it is well known that over an
$\aleph _0$-noetherian ring any module of projective dimension at
most $n$ is filtered by countably generated (presented) modules of
projective dimension at most $n$. Specializing  to the case of
projective dimension at most one, we observe in
Proposition~\ref{P:Q-aleph_0noeth} that for right orders in $\aleph
_0$-noetherian rings $(\PP_1, \PP_1^\perp)$ is of countable type.

 For rings with a two-sided (Ore) classical
ring of quotients $Q$ we look for descent type results. We consider
the problem of getting information on $(\PP_1, \PP_1^\perp)$
assuming that the right $Q$-modules of finite projective dimension are exactly
the projective $Q$-modules.

We recall that, for a general ring $R$, the right small finitistic
dimension, f.dim $R$ is the supremum of the projective dimension of
modules in $\rmod R$ with finite projective dimension. The big finitistic
dimension, F.dim $R$ is the supremum of the projective dimension of
right $R$-modules of finite projective dimension.

In Proposition~\ref{P:divisible} we show that if $R$ has a two-sided
classical ring of quotients $Q$ such that f.dim $Q=0$, then
$\PP_1(\rmod R)^\perp$ coincides with the $\mathrm{Ext}$-orthogonal
of the cyclically presented modules in $\PP_{1}(\rmod R)$.
Therefore, in this case, if $(\PP_1, \PP_1^\perp)$ is of finite
type, then the modules in $\PP_1$ are precisely the direct summands
of modules filtered by the cyclically presented modules of
projective dimension at most one, and $\PP_1^\perp=\mathcal{D}$ the
class of  divisible modules.

To work with a countably presented module $M\in \PP _1$ we use the
relative Mittag-Leffler conditions, that first appeared in \cite{BH}
and were further developed in \cite{AH}, as an effective way to
characterize vanishing conditions of the functor $\mathrm{Ext}$.

In Theorem~\ref{P:1} we patch together the results for countably presented modules
with the ones giving the countable type proving that if the ring is an order in an $\aleph
_0$-noetherian ring $Q$ of small finitistic dimension $0$, then  the
cotorsion pair generated by the modules of projective dimension at
most one is of finite type if and only if $Q$ has big finitistic
dimension $0$.  As a consequence of our
work we find, for example, that $(\PP_1,\mathcal{D})$ is a cotorsion
pair of finite type for orders in semisimple artinian rings
(Corollary~\ref{C:goldie}) so, in particular, for commutative
domains (Corollary~\ref{C:domains}). We also characterize the
commutative noetherian rings for which $(\PP_1,\PP_1^\perp)$ is of
finite type as the ones that are orders into artinian rings.

We remark that this kind of results had been only considered in the commutative domain setting. The cotorsion pair $(\PP_1,\PP_1^\perp)$ was known to be
of finite type only in these two cases: the class of
Pr\"ufer domains and the class of Matlis domains. For the first
class the key result  \cite[VI Theorem 6.5]{FS} is  that a module of
projective dimension at most one over a Pr\"ufer domain is
filtered by cyclic finitely presented
modules (which are all of projective dimension one). For the second
class, recall that a domain $R$ is a Matlis domain provided that the
quotient field $Q$ of $R$ has projective dimension one. If this is
the case, then Matlis proved that the class of divisible module
coincides with the class of epimorphic images of injective modules
(see \cite{M}). From this fact it easily follows that $^\perp
\D=\PP_1$ (see
also \cite{L}).

The paper is structured as follows, in Section~\ref{preliminaries}
we introduce notations and some basic facts about cotorsion pairs.
The notions concerning relative Mittag-Leffler modules   are given
in Section~\ref{ML}, where we also prove the results about these
modules which will be needed in the sequel. We specialize to modules
of bounded projective dimension in Section~\ref{pn}, and
we examine the question of the countable type in
Section~\ref{countable}.

In Sections~\ref{quotients} and \ref{sufficient}, we assume that $R$
has a classical ring of quotients with finitistic dimension $0$ and
we investigate the consequences on the class $\PP _1$, proving
Theorem~\ref{P:1} which is the main result of this part of the
paper. We devote Section~\ref{S:positive} to expose some
applications of our work, and we finish in Section~\ref{examples}
with a discussion of examples and counterexamples that limit the
scope for possible generalizations. In particular, we exhibit
examples showing that $(\PP_n, \PP_n^\perp)$  of finite type does not
imply the finite type of $(\PP_{n-1}, \PP_{n-1}^\perp)$ (Example~\ref{E:Birge} and
Proposition~\ref{P:Birge}).

{\bf Acknowledgments.} We kindly thank Santiago Zarzuela for
providing us with Examples~\ref{E:counterex} (i) and Birge
Huisgen-Zimmermann for suggesting the use of Example~\ref{E:Birge}.

\section{Preliminaries and notations}\label{preliminaries}

Let $R$ be an associative ring with unit.

For any infinite cardinal $\mu$, $\mk{\mu}$ and $\mk{<\mu}$ will be
the classes  of modules with a projective resolution consisting of
$\leq\mu$-generated or $<{}\mu$-generated projective modules,
respectively.  We will simply write $\rmod R$ for $\mk{<\aleph_0}$.
For any class $\C$ of right (left) $R$-modules, $\C(\rmod R)$ and
$\C( \mk{\aleph_0})$ will denote the classes $\C\cap \rmod R$ and
$\C\cap \mk{\aleph_0}$, respectively.

An ascending chain $(M_{\alpha}\mid \alpha<\mu)$ of submodules of a module $M$ indexed by a
cardinal $\mu$ is called \emph{continuous} if $M_{\alpha}=\cup_{\beta <\alpha}M_{\beta}$ for all
limit ordinals $\alpha<\mu$. It is called a filtration of $M$ if $M_0=0$ and $M=\cup_{\alpha
<\mu}M_{\alpha}$.

Given a class $\C$ of modules, we say that a module $M$ is $\C$-\emph{filtered} if it admits a
filtration  $(M_{\alpha}\mid \alpha<\mu)$ such that $M_{\alpha+1}/M_{\alpha}$ is isomorphic to
some module in $\C$ for every $ \alpha<\mu$. In this case we say that $(M_{\alpha}\mid
\alpha<\mu)$ is a $\C$-\emph{filtration} of $M$.

For every class $\C$ of right $R$-modules we set
\[\C^\perp=\{X\in \Mod R\mid \Ext_R^i(C, X)=0 \ {\rm for \ all}\ C\in \C\ {\rm for \ all}\ i\geq 1 \}\]
\[^\perp \C=\{X\in \Mod R\mid \Ext_R^i(X, C)=0 \ {\rm for \ all}\ C\in \C\ {\rm for \ all}\ i\geq 1 \}\]
\[\C^{\perp_1}=\{X\in \Mod R\mid \Ext_R^1(C, X)=0 \ {\rm for \ all}\ C\in \C \}\]
\[^{\perp_1} \C=\{X\in \Mod R\mid \Ext_R^1(X, C)=0 \ {\rm for \ all}\ C\in \C \}\]
A pair of classes of modules $(\A, \B)$ is a \emph{%
cotorsion pair} provided that $\A = {}^{\perp_1}\B$
 and $\B =\A ^{\perp_1}$. Note that for every class $\C$, ${^\perp \C}$ is a \emph{resolving} class, that is,
 it is closed under extensions, kernels of epimorphisms and contains the projective modules. In particular,
 it is syzygy-closed. Dually, $\C ^\perp$ is \emph{coresolving}: it is closed under extensions,
 cokernels of monomorphisms and contains the injective modules. In particular, it is cosyzygy-closed.
 A pair $(\A, \B)$ is called a \emph{hereditary cotorsion pair} if $\A={^\perp \B}$ and $\B=\A ^\perp $. It is easy
 to see that $(\A, \B)$ is a hereditary cotorsion pair if and only if $(\A,\B)$ is a cotorsion pair such that $\A$
 is resolving, if and only if $(\A,\B)$ is a cotorsion pair such that $\B$ is coresolving.

A cotorsion pair $(\A, \B)$ is \emph{complete} provided that every right $R$-module $M$ admits a \emph{special $\A$-precover}, that is, if there exists an exact sequence of the form $0\to B \to A \to M\to 0$  with $B\in \B$ and $A\in \A$.
For a class $\mathcal{C}$ of right modules, the pair
$(^{\perp }(\C^{\perp}),\C^{\perp})$ is a (hereditary) cotorsion pair; it is called the cotorsion pair
\emph{generated} by $\C$. Clearly,
${^{\perp} (\mathcal C ^{\perp})} = {^{\perp_1} (\mathcal C ^{\perp_1})}$ provided that a first syzygy of
$M$ is contained in $\mathcal C$ whenever $M \in \mathcal C$.

Every cotorsion pair generated by a set of modules is complete, \cite{ET}. If all the modules in $\C$ have projective
dimension $\leq n$, then $^{\perp }(\C^{\perp}) \subseteq \PP_n$ as well.

In computing Ext-orthogonals of $\C$-filtered modules the following,
known as Eklof's Lemma,  is essential.

\begin{Fc} \cite[XII.1.5]{EM}\label{F:Eklof-filt}
Let $R$ be a ring and $M$, $N$ be right $R$-modules. Assume that $M$
has a
filtration $(M_{\alpha }\mid \alpha <\sigma )$ such that $\hbox{Ext}%
_{R}^{1}(M_{\alpha +1}/M_{\alpha },N)=0$ for all $\alpha +1<\sigma $. Then
$\hbox{Ext}_{R}^{1}(M,N)=0$.
\end{Fc}

We recall also the following useful description of the modules in the first component of a cotorsion pair

\begin{Fc} \cite[Theorem 2.2]{T}\label{F:Jan-description} Let $\C$ be a set of right $R$-modules.  An $R$-module belongs
to the cotorsion pair generated by $\C$ if and only if it is a direct summand of a $\mathcal C'$-filtered module
where $\mathcal C' = \mathcal C \cup \{R\}$.
\end{Fc}

A hereditary cotorsion pair $ (\A,\B)$ in $\Mod R$ is of \emph{countable type} ( \emph{finite type})
provided that there is a class $\ES $ of modules in $\mk{\aleph_0} $ ($\rmod R$)
such that $\ES$ generates $(\A, \B)$, that is $\ES ^\perp= \B$.

We denote by $\PP$ the class of right $R$-modules of finite
projective dimension,  and for every $n\geq 0$, we denote by $\PP_n$
the class of right $R$-modules of projective dimension at most $n$.
In case we need to stress the ring $R$ we shall write $\PP (R)$ and
$\PP_n (R)$, respectively.

In \cite{AEJO} it is shown that, for every $n\in \N$,  $(\PP_n, \PP_n^\perp)$ is a hereditary
cotorsion pair; moreover it is complete, since it is generated by a set of representatives of the
modules in the  class $\PP_n(\mk{\mu})$ where $\mu=\max \{\card R, \aleph_0\}$.

We consider also the cotorsion pair generated by the class
$\PP_n(\rmod R)$, for every $n\in \N$ that is the cotorsion pair
$(^\perp{}(\PP_n(\rmod R)^\perp), \PP_n(\rmod R)^\perp) $. By
definition, this  cotorsion pair is of finite type; it is also
hereditary because the class $\PP_n(\rmod R)$ is resolving. Clearly,
the class $^\perp{}(\PP_n(\rmod R)^\perp)$ is contained in $\PP_n$.

We are interested in $n$-tilting cotorsion pairs and in cotorsion pairs associated to subclasses of $\PP_n$.

Recall that an $n$-tilting cotorsion pair is the hereditary cotorsion pair $(^\perp (T^\perp),
T^\perp)$ generated by an $n$-tilting module $T$. The class $T^\perp$ is then called $n$-tilting
class.  If $\ES$ is a subclass of $\PP_n(\rmod R)$ then the hereditary cotorsion pair generated by $\ES$, that is $(^{\perp }(\ES^{\perp}),\ES^{\perp})$, is an $n$-tilting cotorsion pair. By results in \cite{BET}, \cite{BH} , \cite{ST1} and \cite{BS} all $n$-tilting cotorsion pairs can be generated in this way, namely they are of finite type. Consequently, the
class $(\PP_n(\rmod R))^\perp$ is the smallest $n$-tilting class.

We will consider also $\Tor$ orthogonal classes.
For every class $\C$ of right $R$ modules we set
\[\C^{\intercal}=\{X\in R \LMod \mid \Tor_i^R(C, X)=0 \ {\rm for \ all}\ C\in \C\ {\rm for \ all}\ i\geq 1 \}\]

\section{Relative Mittag-Leffler conditions} \label{ML}

The definition of Mittag-Leffler inverse systems goes back to Grothendieck \cite[Proposition 13.1.1]{Gr}. Raynaud and Gruson in \cite{RG}  realized  the strong connection between this concept and the notion of Mittag-Leffler modules.

We recall here a weaker notion, that is the Mittag-Leffler condition restricted to particular classes.

\begin{Def}
 Let $M$ be a right module over a ring $R$, and let
$\Q$ be a class of left $R$-modules. We say that $M$ is a
\emph{$\Q$-Mittag-Leffler module} if the canonical map
\[\rho \colon  M\bigotimes_R \prod _{i\in I}Q_i\to \prod _{i\in
I}(M\bigotimes_RQ_i)\] is injective for any family $\{Q_i\}_{i\in I}$ of modules in $\Q$.
\end{Def}

Taking $\Q=R \lMod $ we recover Raynaud and Gruson's notion of Mittag-Leffler modules.

We shall use the following characterization of relative
Mittag-Leffler modules.

\begin{Th} \label{AF2} {\rm (\cite[Theorem 5.1]{AH})} Let $\Q$ be a class of left $R$-modules. For a right
$R$-module $M$,   the following statements are equivalent:
\begin{itemize}
\item[(1)] $M$   is  $\Q$-Mittag-Leffler.
\item[(2)] Every
direct system of  finitely presented right $R$-modules
$(F_\alpha,u_{\beta \,\alpha})_{\beta \,\alpha \in I  }$    with
$M=\varinjlim (F_\alpha,u_{\beta \alpha})_{\beta , \alpha \in I }$
has the property that for  any $\alpha \in I  $ there exists $\beta
\ge \alpha$ such that, for any $Q\in \Q$, $\mathrm{Ker}\,(u_{\beta
\alpha}\otimes _RQ)=\mathrm{Ker}\,(u_{\alpha}\otimes _R Q)$, where
$u_\alpha \colon  F_\alpha \to M$ denotes the canonical map.
\item[(3)] There exists
direct system of  finitely presented right $R$-modules
$(F_\alpha,u_{\beta \,\alpha})_{\beta \,\alpha \in I  }$    with
$M=\varinjlim (F_\alpha,u_{\beta \alpha})_{\beta , \alpha \in I }$
satisfying that for  any $\alpha \in I  $ there exists $\beta \ge
\alpha$ such that, for any $Q\in \Q$, $\mathrm{Ker}\,(u_{\beta
\alpha}\otimes _RQ)=\mathrm{Ker}\,(u_{\alpha}\otimes _R Q)$, where
$u_\alpha \colon F_\alpha \to M$ denotes the canonical map.
\end{itemize}
\end{Th}

The relation between relative Mittag-Leffler modules and cotorsion
pairs of finite type is given by the following result.
\begin{Th}{\rm ( \cite[Theorem 9.5]{AH}, \cite[Theorem~5.1]{BH})} \label{T:Dol-Lid} Let $R$ be an arbitrary ring. Let $(\A, \B)$ be a cotorsion pair of
finite type. Let $\ES=\A\bigcap \rmod R$ and let $\C=\A^{\intercal}.$ Then
every module in $\A$ is $\C$-Mittag-Leffler. If $M$ is   a countably
presented right $R$-module that is a direct limit of modules in
$\ES$, then $M$ is in $\A$ if and only if it is $\C$-Mittag-Leffler.
\end{Th}

We illustrate now some properties of $\Q$-Mittag-Leffler modules
that will be used later on. First we state an auxiliary Lemma.

\begin{Lemma}\label{simple} Let $\mu \colon A\to B$ a morphism of right $R$-modules. Let $A'$ and $B'$ denote submodules
of $A$ and $B$, respectively, such that $\mu (A')\subseteq B'$. Let
$\mu '\colon A'\to B'$ be the restriction of $\mu$. Then the kernel
of the induced map $f\colon B'/\mu '(A')\to B/\mu (A)$ is
$\mathrm{ker}\, f=\left( \mu(A)\cap B'\right)/\mu'(A')$.\end{Lemma}
\begin{Prop}\label{P:submodules} Let $R$ be a ring, and let $M_R\in \PP_1$. Assume that $M$ is a $\Q$-Mittag-Leffler
module where $\Q$ is a class of left $R$-modules contained in $ M^\intercal$.
Let $\Y$ be a class of left $R$-modules consisting of submodules of modules in $\Q$ such
that $\Y\subseteq M^\intercal$. Then $M$ is a $\Y$-Mittag-Leffler module.
\end{Prop}
\begin{Proof} Using the Eilenberg trick, if needed, we can
assume that $M$ has a free presentation
\[(1)\quad 0\to R^{(J)}\overset{\mu}\to R^{(I)}\to M\to 0,\]
where   $I$ and $J$ are sets.

Since finitely presented modules are Mittag-Leffler, we can assume
that either $I$ or $J$ is infinite. As we are stating a property
on $M^\intercal$ and free modules are Mittag-Leffler, we can cancel
the free direct summands of $M$. Hence, without loss of generality,
we may assume that the image of $\mu$ has non zero projection on all
the direct summands of $R^{(I)}$, and therefore that $J$ is an
infinite set.

For every finite subset $F$ of $J$, let $\mu_{F}$ be the restriction
of $\mu$ to $R^F$ and let $G_F$ be the smallest subset of $I$ such
that $\mu _F(R^F)\leq R^{G_F}$. Let $C_F$ be the finitely presented
right $R$-module $R^{G_F}/\mu_F(R^F)$; then $C_F\in \PP_1(\rmod R)$.
Let $\F$ be the family of the finite subsets of $J$ and consider the
direct system $(C_F; f_{KF})_{F\subseteq K\in \F}$ where the
structural morphisms $ f_{KF} \colon C_F\to C_K$ are induced by the
injections of $R^{G_F}$ into $R^{G_K}$. Then, $M_R$ is isomorphic to
the direct limit of the direct system  $(C_F; f_{KF})_{F\subseteq
K\in \F}$. Let $f_F\colon C_F\to M\cong \varinjlim_FC_F$ be the
canonical morphisms.

For every $F\leq K\in \F$ and every left $R$-module $N$, we have a commutative
diagram:

\[ \xymatrix{C_F\otimes_R N\quad\ar[r]^{\small{f_F\otimes_R 1_N}} \ar[d]_{\small{f_{KF}\otimes_R 1_N} }
                               &M\otimes_R N\\
                                 C_K\otimes_R N\ar[ur]_{\small{f_K\otimes_R 1_N}}}.\]

By the definitions of the finitely presented modules $C_F$ and of
the maps $f_F$ and $f_{KF}$, Lemma~\ref{simple} allows us to
conclude
\[(a)\qquad\ker (f_F\otimes_R 1_N)=\dfrac{\mu\otimes_R 1_N \left (N^{(J)}\right)\cap N^{G_F}}{\mu_F\otimes_R 1_N \left (N^F\right )}\]
and
\[(b)\qquad \ker (f_{KF}\otimes_R 1_N)=\dfrac{\mu_K\otimes_R 1_N \left (N^K\right)\cap N^{G_F}}{\mu_F\otimes_R 1_N \left (N^F\right )}\]
where, for any set $L$, we identify $R^{(L)}\otimes _RN$ with
$N^{(L)}$.

By Theorem~\ref{AF2}(2), the assumption that $M$ is a
$\Q$-Mittag-Leffler module amounts to the following

(*) for every, $F\in \F$ there is a subset $l(F)\in \F$, $l(F)\geq F$ such that
\[\left[\mu\otimes_R 1_Q \left (Q^{(J)}\right )\right]\cap Q^{G_F}=
\left[\mu_{l(F)}\otimes_R 1_Q \left (Q^{l(F)}\right )\right]\cap
Q^{G_{F}},\] for every $Q\in \Q$.

 Let now  ${}_RY\in \Y$ be a submodule of some module $Q\in \Q$. We claim that
\[\left[\mu\otimes_R 1_Y \left (Y^{(J)}\right )\right]\cap Y^{G_F}=
\left[\mu_{l(F)}\otimes_R 1_Y \left (Y^{l(F)}\right )\right]\cap
Y^{G_{F}}.\]

Observe that only the inclusion $\subseteq $ of the claim needs to
be proved. Consider the commutative diagram:
\[ \begin{CD}
@. 0@. 0\\
@. @VVV @VVV\\
0@>>>Y^{(J)}@>\small{\mu\otimes_R 1_Y}>>Y^{(I)} @>>>M\otimes_R Y@>>> 0\\
@. @V{\sigma}VV @VV{\tau}V@VVV \\
0 @>>>Q^{(J)}@>>\small{\mu\otimes_R 1_{Q}}> Q^{(I)} @>>>
M\otimes_RQ@>>> 0
 \end{CD}\]
 where the rows are exact by the hypothesis that $Q, Y\in M^\intercal$

Condition (*) and the commutativity of the above diagram yield:
\[\tau \left(\mu\otimes_R 1_Y\left(Y^{(J)}\right)\cap Y^{G_F} \right)=\tau
\left (\mu\otimes_R 1_Y\left (Y^{(J)}\right) \right)\cap \tau(Y^{G_F})\leq \mu\otimes_R 1_Q \left (Q^{(J)}\right )
\cap\  Q^{G_F}=\]
\[= \mu_{l(F)}\otimes_R1_Q \left (Q^{l(F)}\right )\cap Q^{G_{F}}\]

Let $\underline y\in Y^{(J)}$ be such that $\mu\otimes_R
1_Y(\underline y)\in Y^{G_F}$. By the above inclusion,
\[\tau \left (\mu\otimes_R
1_Y(\underline y)\right)=\mu_{l(F)}\otimes_R1_Q(\underline z),\] for
some $\underline z\in Q^{l(F)}$ with
$\mu_{l(F)}\otimes_R1_Q(\underline z)\in Q^{G_{F}}$. Since
$\mu_{l(F)}\otimes_R 1_Q$ is the restriction of $\mu\otimes_R1_Q$,
 $\mu_{l(F)}\otimes_R1_Q(\underline
z)=\mu\otimes_R1_Q(\underline z)$. Thus,
\[\tau\left (\mu\otimes_R1_Y(\underline y)\right)=\mu\otimes_R1_Q\left(\sigma (\underline y)\right)=
\mu\otimes_R1_Q(\underline z)\in Q^{G_{F}} .\] By the injectivity of
$\mu\otimes_R1_Q$ we conclude that $\sigma(\underline y)=\underline
z$, hence $\underline y\in Y^{l(F)}$.  This proves the claim.

Then, taking into account (a) and (b), we conclude that $M$ is a
$\Y$-Mittag-Leffler module using Theorem~\ref{AF2}(3).
\end{Proof}

Before giving other properties of relative Mittag-Leffler modules,
we need a lemma.

\begin{Lemma}\label{L:presentation} Let $M$ be a right $R$-module and   let $\mu$ be an infinite
cardinal. $M$ is $<\mu$-presented if and only if there exists a direct system $(C_\alpha,
u_{\beta\, \alpha}\colon C_\alpha \to C_\beta)_{\alpha \le \beta \in \Lambda}$ such that
$M=\varinjlim C_\alpha $ and the cardinality of $\Lambda$ is strictly smaller than $\mu$.
\end{Lemma}

\begin{Proof} If $\mu=\aleph_0$ the claim is obvious. Assume that $M$ is $<\mu$-presented and $ \mu>\aleph_0$. Let $\{x_i\}_{i\in I}$ be a generating set of
$M$ such that $\abs I <\mu$. Consider the exact sequence
\[0\to L\stackrel{g}{\to}R^{(I)}\stackrel{f}{\to}M\to 0\]
where, if $\{e_i\}_{i\in I}$ denotes the canonical basis of $R^{(I)}$, $f(e_i)=x_i$ for any $i\in
I$. By hypothesis we can choose a generating set  $\{y_j\}_{j\in J}$ of $L$ such that $\abs{J} <\mu$.

For any finite subset $F$ of $J$ there exists a finite subset $I(F)$ of $I$ such that $g(\sum
_{j\in F}y_jR)\subseteq R^{I(F)}$. Setting $C_F= R^{I(F)}/\sum _{j\in
F}y_jR$, we obtain a direct system of finitely presented modules with limit $M$ indexed by the set
$\F$ of finite subsets of $I$. $\F$ has less than $\mu$ elements.

For the converse,  let $(C_\alpha, u_{\beta\, \alpha})_{\alpha \le \beta \in \Lambda}$ be a direct
system of finitely presented modules such that $M=\varinjlim C_\alpha $; assume $\abs \Lambda
<\mu$. The canonical presentation of the direct limit (see \cite[Proposition 2.6.8]{weibel})
\[\bigoplus _{\alpha \le \beta} C_{{\beta\, \alpha}}\stackrel{\Phi}{\to}\bigoplus _{\alpha \in \Lambda}
C_\alpha \to M \to 0\]
where for every $\alpha \le \beta$, $C_{{\beta\, \alpha}}= C_\alpha$, gives a pure exact sequence
\[0\to \mathrm{Im} \Phi \to\bigoplus _{\alpha \in \Lambda}
C_\alpha \to M \to 0. \]

Since $\oplus _{\alpha \le \beta} C_{{\beta\, \alpha}}$ is $<\mu$-generated, so is $\mathrm{Im}
\Phi $. Moreover, since  $\oplus _{\alpha \in \Lambda} C_\alpha $ is $<\mu$-presented we conclude that
$M$ is $<\mu$-presented.
\end{Proof}

\begin{Prop}\label{P:coherence} Let $\mu$ be an infinite cardinal, and let
$M$ be a $<\mu$-generated $R$-Mittag-Leffler right $R$-module. Then $M$ is $<\mu$-presented.
\end{Prop}

\begin{Proof}
Assume first that $\mu =\aleph _0$, so that $M$ is a finitely generated module. To show
that $M$ is finitely presented we only need to show that the natural map
\[\rho _J\colon M\otimes R^J\to (M\otimes R)^J\]
is bijective for any set $J$ (cf.\cite[Theorem 3.2.22]{enochsjenda}). Since $M$ is finitely generated, for any
set $J$, $\rho _J$ is onto \cite[Lemma 3.2.21]{enochsjenda} and by our assumption  $\rho _J$ is
injective, hence bijective.

Assume now $\mu >\aleph _0$.  Let $X=\{x_i\}_{i\in I}$ be a set of generators of $M$. As $M$ is $R$-Mittag-Leffler,
for any finite subset $F$ of $I$ there exists a submodule $N_F$ of $M$ that is countably presented
and contains $\{x_i\}_{i\in F}$ (cf. \cite[Theorem 5.1~(4)]{AH}).

Let $\mathcal{F}$ denote the set of all finite subsets of $I$, then $M$ is the directed union of
$(N_F)_{F\in \mathcal{F}}$. As each $N_F$ is a countable direct limit of finitely presented
modules, we deduce that $M$ is the direct limit of a direct system of
finitely presented modules indexed by a set of the same cardinality as $\mathcal{F}$. Since $\mathcal{F}$ has cardinality $<\mu$, this implies by
Lemma~\ref{L:presentation} that $M$ is $<\mu$-presented.
\end{Proof}

\begin{Cor}\label{C:coherence} Let $(\A ,\B)$ be a hereditary cotorsion pair of finite type,
and let $\mu$ be an infinite cardinal. If $M$ is a right $R$-module of $\A$ that is
$<\mu$-generated then $M\in \A (\mk{<\mu})$
\end{Cor}

\begin{Proof} First  observe that, since the cotorsion pair, is hereditary the class $\A$ is
resolving; so to prove the statement it is enough to show that if $M\in \A$ is $<\mu$-generated
then it is $<\mu$-presented.

Since $M\in \A$, it is $\A
^\intercal$-Mittag-Leffler  by Theorem~\ref{T:Dol-Lid}. Hence $M$ is  $R$-Mittag-Leffler and thus the conclusion follows by Proposition~\ref{P:coherence}.
\end{Proof}

By \cite[Proposition~9.2]{AH} the conclusion of
Corollary~\ref{C:coherence} holds, more generally, for hereditary
cotorsion pairs $(\A ,\B)$ of countable type and such that the class
$\B$ is closed by direct sums.
\section{The cotorsion pair $(\PP_n, \PP_n^\perp)$} \label{pn}

We start characterizing when this cotorsion pair is of finite type.

\begin{Prop}\label{P:finite} Let $R$ be a ring. The following conditions are equivalent:
\begin{enumerate}
\item[(i)] The class $\PP_n^\perp$ is closed under direct sums.
\item[(ii)]  $(\PP_n, \PP_n^\perp)$ is an $n$-tilting cotorsion pair.
\item[(iii)] The cotorsion pair $(\PP_n, \PP_n^\perp)$ is of finite type.
\item[(iv)] $\PP_n^\perp=\PP_n(\rmod R)^\perp$.
\item[(v)] Every module in $\PP_n$ is a direct summand of a $\PP_n(\rmod R)$-filtered module.
\end{enumerate}
When the above equivalent conditions hold, then  $\PP _n\subseteq
\varinjlim \PP_n(\rmod R)$.
\end{Prop}

\begin{Proof} (i) $\Rightarrow$ (ii). A hereditary cotorsion pair $(\A, \B)$ is an $n$-tilting cotorsion pair if and only if it is complete, $\A\subseteq \PP_n$ and $\B$ is
closed under direct sums (see \cite{AC},  \cite{KS} or \cite{T2}).  Since $(\PP_n, \PP_n^\perp)$ is a complete cotorsion pair, condition (i) implies (ii).

(ii) $\Rightarrow$ (iii). Any $n$-tilting cotorsion pair $(\A, \B)$ is of finite type, by  \cite{BH} and \cite{BS}.

(iii) $\Leftrightarrow$ (iv). By definition, a cotorsion pair $(\A,
\B)$ is of finite type if and only if it is generated by
(representatives of) the modules in $\A(\rmod R)$.

(iv)$\Leftrightarrow$ (v). Is a consequence  of
Fact~\ref{F:Jan-description}.

(iii) $\Rightarrow$ (i). This follows by the fact that  for every $M\in \rmod R$, the functors
$\Ext^i_R(M, -)$ commutes with direct sums.

If the conditions hold, then the rest of the claim follows from
\cite[proof of Theorem~2.3]{AT2}.
\end{Proof}



Trivially, $(\PP _0, \PP _0 ^\perp )$ is of finite type. Note that,
in this case, condition (v) in Proposition~\ref{P:finite} can be
stated by saying that any projective right module is a direct
summand of an $R$-filtered (hence free) module.

It is well known that $\PP_1\subseteq  \varinjlim \PP_1(\rmod R)$.
This can be seen as a consequence of the fact that $(\PP _0, \PP
_0^\perp)$ is of finite type. The rest of this section will be
devoted to extend this result to arbitrary projective dimension.
That is,   $(\PP_{n-1}, \PP_{n-1}^\perp)$   of finite type implies
$\PP_n\subseteq \varinjlim \PP_n(\rmod R)$. Our arguments follow the
ones  in \cite{BS}.

\medskip
First we state a Lemma.

\begin{Lemma}\label{L:lim-mod} Let $R$ be a ring. Let
\[0\to H\to G\to C\to 0\]
be an exact sequence of right $R$-modules. Let $\mu$ be an infinite
cardinal. Then,
\begin{itemize}
\item[(i)] if there exists $n\ge 0$ such that  $H$ and $C$ are in $\PP_{n}(\mk{<\mu})$ then also
$G\in \PP_{n}(\mk{<\mu})$.

\item[(ii)] If $H$ and $G$ in $\PP_{n-1}(\mk{<\mu})$, for some $n\geq 1$, then $C\in \PP_{n}(\mk{<\mu})$.
\end{itemize}
\end{Lemma}
\begin{Proof} Statement $(i)$ follows by inductively applying the
Horseshoe Lemma.

To prove (ii) we can assume that $n>1$. Let $0\to G_1\to P_0\to G\to
0$ be an exact sequence with $P_0$ a $<\mu$-generated projective
module and $G_1\in \PP_{n-2}(\mk{<\mu})$. By a pull-back  argument
we obtain the following commutative diagram:
\[\begin{CD}
@. 0@. 0\\
@. @AAA @AAA\\
0@>>>H@>>>G@>>>C@>>> 0\\
@. @AAA @AAA@|\\
0@>>>X@>>>P_0@>>>C@>>> 0\\
@. @AAA @AAA@. \\
@. G_1@=G_1 @.\\
@. @AAA @AAA\\
@. 0@. 0
 \end{CD}.\]
Applying $(i)$ to the exact sequence $0\to G_1\to X\to H \to 0$ we
deduce that  $X\in \PP_{n-1}(\mk{<\mu})$. Hence  $C\in
\PP_{n}(\mk{<\mu})$.
\end{Proof}

Following the ideas in  \cite{BS}, we look at conditions on the
syzygy module of $M\in \PP_n$. To this aim, we state a result for
$\C$-filtered modules, where $\C$ is a class of $<\mu$-presented
modules for some infinite cardinal $\mu$. The proof of this result
for the case of $\mu\geq \aleph_1$ appears in \cite[XII.1.14]{EM}
and in \cite[Proposition 3.1]{BS} for the case $\mu=\aleph_0$. An
alternative proof is in \cite[Theorem 6]{ST}.

\begin{Prop}(\cite[XII 1.14]{EM}, \cite[Prop. 3.1]{BS}, \cite[Theorem 6]{ST}) \label{P:Eklof} Let $\mu$ be an infinite
cardinal. Let $M$ be a $\C$-filtered  module where $\C$ is a family of $<\mu$-presented modules. Then there exists a subset $\ES$ of $\C$-filtered submodules of $M$ satisfying the following properties:
\begin{enumerate}
\item[(1)] $0\in \ES$.
\item[(2)]  $\ES$ is closed under unions of arbitrary chains.
\item[(3)]  For every $N \in \ES$, $N$ and $M/N$ are $\C$-filtered.
\item[(4)]  For every subset $X \subseteq M$ of cardinality $<\mu$, there is a $<\mu$-presented module $N \in \ES$ such that $X \subseteq N$.
\end{enumerate}
\end{Prop}
An immediate consequence of conditions (2) and (4) in Proposition~\ref{P:Eklof} is the following.

\begin{Cor} \label{C:Eklof-} Let $\mu$ be an infinite
cardinal. Let $M$ be a $\mu$-generated $\C$-filtered  module where
$\C$ is a family of $<\mu$-presented modules. Then there is a
filtration $(M_{\alpha}\mid \alpha\in \mu)$ of $M$ consisting of
$<\mu$-presented submodules of $M$ such that $M_{\alpha}$ and
$M/M_{\alpha}$ are $\C$-filtered for every $\alpha\in \mu$.
\end{Cor}

The next result is a straight generalization of \cite[Lemma 3.3]{BS}

\begin{Lemma}\label{L:filtration} Let $\mu$ be an infinite cardinal, and let $\C$ be a family of $<\mu$-presented right modules
containing the regular module $R$. Let $M$ be a $\mu$-presented
right module, and let \[0\to K \to F\to M\to 0\] be a free
presentation of $M$ with $F$ and $K$ $\mu$-generated. Assume that
$K$ is  a direct summand of a $\C$-filtered  module. Then, there
exists an exact sequence:
\[0\to H\to G \to M\to 0\]

where $H$ and $G$ are $\mu$-generated $\C$-filtered modules.

\end{Lemma}

\begin{Proof} Let $K$ be a summand of a $\C$-filtered module $P$. Since $K$ is $\mu$-generated, Proposition~\ref{P:Eklof} implies that $K$ is contained in a $\mu$-generated $\C$-filtered submodule of $P$; thus we may assume that $P$ is $\mu$-generated. By Eilenberg's trick, $K\oplus P^{(\omega)}\cong P^{(\omega)}$. Consider the exact sequence
\[0\to K\oplus P^{(\omega)}\to F\oplus P^{(\omega)}\to M\to 0\]
and let $H=K\oplus P^{(\omega)}\cong  P^{(\omega)}$, $G=F\oplus P^{(\omega)}$. Then $G$ and $H$
are $\mu$-generated $\C$-filtered modules.
\end{Proof}

Now we are ready to  prove the announced result.

\begin{Th}\label{T:lim-in P_n} Let $R$ be a   ring, and let $n\ge 1$.
\begin{enumerate}
\item[(i)] If the cotorsion pair generated by $\PP_{n-1}(\mk{\aleph_0})$ is of finite type, then every module in $\PP_{n}(\mk{\aleph_0})$ is a direct limit of modules in $\PP_n(\rmod R)$.
\item[(ii)] If the cotorsion pair
$(\PP_{n-1}, \PP_{n-1}^\perp)$ is of finite type, then every module in $\PP_n$ is a direct limit
of modules in $\PP_n(\rmod R)$.
\end{enumerate}
 \end{Th}
\begin{Proof} Statements $(i)$ and $(ii)$ are clear for $n=1$. Hence we may assume that $n>1$.

(i) Let $M\in \PP_{n}(\mk{\aleph_0})$.  Then there is an exact sequence
\[ 0\to K\to F_0\to M\to 0\]
where $F_0$ is an $\aleph_0$-generated free module and   $K\in \PP_{n-1}(\mk{\aleph_0})$. By
assumption $K$ is a direct summand of a $\PP_{n-1}(\rmod R)$-filtered module.

By Lemma~\ref{L:filtration} applied to the family $\PP_{n-1}(\rmod R)$ for the case
$\mu=\aleph_0$, there exists an exact sequence
\[ 0\to H\to G\to M\to 0\]
where $H$ and $G$ are countably generated $\PP_{n-1}(\rmod R)$-filtered modules. By
Corollary~\ref{C:Eklof-}, $H$ and $G$ admit filtrations $(H_i\mid i\in \N)$ and $(G_j\mid j\in
\N)$, respectively, consisting of finitely presented $\PP_{n-1}(\rmod R)$-filtered submodules.
Without loss of generality we can assume that $H$ is a submodule of $G$. Given $i < \omega$, there
is an $j(i)$ such that $H_i\subseteq G_{j(i)}$; and we can choose the sequence $(j(i))_{i<\omega}$
to be strictly increasing. Consider the exact sequence
\[0\to H_i\to G_{j(i)}\to C_i\to 0\]

For every $i\in \N$, the  modules  $H_i$ and $G_{j(i)}$ are finitely presented and they belong to ${}^\perp(\PP_{n-1}(\rmod R)^\perp)$, by Fact~\ref{F:Jan-description}. By Corollary~\ref{C:coherence} they belong to
$\PP_{n-1}(\rmod R)$. Thus, by Lemma~\ref{L:lim-mod}, $C_i\in
\PP_{n}(\rmod R)$. Moreover,
 $M\cong\varinjlim C_i$ by construction, hence (i) follows.

(ii)  By way of contradiction, assume that the result is not true and let $\mu$  be the least
cardinal for which there exists an $R$-module $M\in \PP_{n}(\mk{\mu})$ which is not a direct limit
of modules in $\PP_n(\rmod R)$. By (i), $\mu>\aleph_0$.

There exists an exact sequence
 \[0\to K\to F_0\to M\to 0\]
where $F_0$ is a $\mu$-generated free module and   $K\in \PP_{n-1}(\mk{\mu})$. By assumption $K$
is a direct summand of a $\PP_{n-1}(\rmod R)$-filtered module.

By Lemma~\ref{L:filtration} applied  to the family $\PP_1(\rmod R)$, there exists an exact
sequence
\[ 0\to H\to G\to M\to 0\]
where $H$ and $G$ are $\mu$-generated $\PP_{n-1}(\rmod R)$-filtered modules. By
Corollary~\ref{C:Eklof-}, $H$ and $G$ admit filtrations $(H_{\alpha}\mid \alpha \in \mu)$ and
$(G_{\alpha}\mid \alpha \in \mu)$, respectively, consisting of $<\mu$-presented $\PP_{n-1}(\rmod
R)$-filtered submodules. Without loss of generality we can assume that $H$ is a submodule of $G$.
Given $\alpha \in \mu$, there is a $\beta(\alpha)\in \mu$ such that $H_{\alpha}\subseteq
G_{\beta(\alpha)}$; and we can choose the sequence $\beta(\alpha)\in \mu$ to be strictly
increasing. Consider the exact sequence
\[0\to H_{\alpha}\to G_{\beta(\alpha)}\to C_{\alpha}\to 0\]

Now, for every $\alpha<\mu$,  the modules $H_{\alpha}$ and
$G_{\beta(\alpha)}$  are $<\mu$-presented and in
${}^\perp(\PP_{n-1}(\rmod R)^\perp)$, by
Fact~\ref{F:Jan-description}. By Corollary~\ref{C:coherence}, they
belong to $\PP_{n-1}(\mk{<\mu} )$. Thus, by Lemma~\ref{L:lim-mod},
$C_{\alpha}\in \PP_n(\mk{<\mu} )$. By the minimality of $\mu$,
$C_{\alpha}$ is a direct limit of objects in $\PP_n(\rmod R)$. Now,
$M\cong\varinjlim_{\alpha\in \mu} C_{\alpha}$, by construction,
hence $M$ is  a direct limit of objects in  $\PP_n(\rmod R)$, too. A
contradiction.
\end{Proof}

\begin{Remark}\label{R:asc-desc} \emph{As $(\PP_{0},
\PP_{0}^\perp)$ is always of finite type, it is easy to find
examples showing  that,  in general, the finite type of $(\PP_{n-1},
\PP_{n-1}^\perp)$ does not imply the finite type of $(\PP_n,
\PP_n^\perp)$. More involved examples will be given in
Examples~\ref{E:counterex}.}

\emph{Moreover, the finite type has not a descent property. In fact,
we will show in Proposition~\ref{P:Birge}, that there exist  artin
algebras with the property that $(\PP_2, \PP_2^\perp)$ is of finite
type, while $(\PP_1, \PP_1^\perp)$ is not.}
\end{Remark}

\section{Countable Type}\label{countable}

We are interested in finding conditions under which the cotorsion pair $(\PP_n, \PP_n^\perp)$ is of finite type.
A necessary condition is that $(\PP_n, \PP_n^\perp)$ be of countable type.  To this regard we recall the following results.

\begin{Fc}\label{F:countable-type} \emph{ If $R$ is a commutative domain, then in \cite[VI 6]{FS} it is proved that every module in $\PP_1$ admits a filtration consisting of countably generated submodules of projective dimension at most one. Hence the cotorsion pair $(\PP_1, \PP_1^\perp)$ is of countable type.}

 \emph{ If $R$ is a right  $\aleph_0$-noetherian ring (that is all the right ideals of $R$ are at most $\aleph_0$-generated), then Raynaud Gruson in \cite[Corollary 3.2.5]{RG} proved that the cotorsion pair  $(\PP_n, \PP_n^\perp)$ is of countable type.  This result appears also in \cite{AEJO} and  \cite[Proposition 2.1]{HS}. }
\end{Fc}

In the one dimensional case,   these two cases can be seen in a
common setting.


\begin{Def} Let $R$ be a ring and let $\Sigma$ denote the multiplicative set of the non zero divisors of $R$.
 A right  $R$-module $D$ is said to be {\sl divisible} if $\Ext^1_R(R/rR, D)=0$, for every  element $r\in \Sigma$.
A left $R$-module $Y$ is said to be {\sl torsion free} if $\Tor^R_1(R/rR, Y)=0$, for every element $r\in \Sigma$. \\
Divisible left $R$-modules and torsion free right $R$-modules are defined in an analogous way.

We denote by $\D$ the class of all divisible right $R$-module and by $\TFcal$ the class of all torsion free left $R$-modules.
\end{Def}

Thus,  a right (left) $R$-module $D$ is divisible if and only if the right (left) multiplication by an element of $\Sigma$ is a surjective map and a left (right) $R$-module $Y$ is torsion free  if and only if the left (right)  multiplication by an element of $\Sigma$ is an injective map.

Moreover, if $\C=\{R/rR\mid r \in \Sigma\}\ \cup\  \{R\}$, then $\D=\C^{\perp}$ and $\TFcal=\C^{\intercal}$.

Examples of torsion free $R$-modules are the submodules of free $R$-modules. If $S$ is a multiplicative subset of $\Sigma$ that satisfies the left  Ore condition, then $S^{-1}R/ R$ is a direct limit of $R/sR$, for $s\in S$. Dually, if $S$ is a multiplicative subset of $\Sigma$ that satisfies the right Ore condition, then $RS^{-1}/ R$ is a direct limit of $R/Rs$ for $s\in S$. Hence we have the following well known fact.

\begin{Lemma}\label{L:torsionfree}  Let $R$ be a ring and let $S$ be a multiplicative subset of $\Sigma$.
\begin{enumerate}
\item[(i)] If $S$  satisfies the left Ore condition, then $\Tor_1^R( S^{-1}R/ R, K)=0$, for any torsion free left $R$-module $K$. In particular, $K$ is embedded in $S^{-1}R\otimes_R K$ via the assignment $y\mapsto 1\otimes_R y$, for any $y\in K$.
\item[(ii)] If $S$  satisfies the right Ore condition, then $\Tor_1^R( K, RS^{-1}/ R)=0$, for any torsion free right $R$-module $K$. In particular,  $K$ is embedded in $K\otimes_RRS^{-1}$ via the assignment $y\mapsto y\otimes_R 1$, for any $y\in K$.

\end{enumerate}
\end{Lemma}

\begin{Lemma}\label{L:countably-generated} Let $R$ be a ring and let $S\subseteq \Sigma$ satisfy the right Ore condition such that $Q=RS^{-1}$ is right $\aleph_0$-noetherian.  Let $F$ be a free right $R$-module  and let $K$ be a submodule of $F$ such that $K\otimes_RQ$ is countably generated as a right $Q$-module. Then, $K$ is contained in a countably generated direct summand of $F$.
\end{Lemma}
\begin{Proof} Let $(e_i; i\in I)$ be a  basis of $F$. For every $i\in I$ denote by $\pi_i\colon F\to e_iR$ the canonical projection. For every subset $X$ of $F$, define the support of $X$ as
 \[\mathrm{supp} (X)=\{i\in I\mid \pi_i(x)\neq 0, \mbox{ for some } x\in X\}.\]
 Choose a set of generators of $K\otimes_R Q$ of the form $\{y_n\otimes_R1\mid n\in \N\}$, where $y_n\in K$ for every $n\in \N$. We claim that $\mathrm{supp} (K) =\mathrm{supp} (\sum_{n\in \N}y_nR)$, hence countable. It is clear that $\mathrm{supp} (\sum_{n\in \N}y_nR)\subseteq \mathrm{supp} (K) $. For the converse, let $y\in K$.
 There exist $r_1, \dots , r_\ell\in R$ and $s\in S$ such that $y\otimes_R1=\sum_{i=1}^{\ell} y_ir_i\otimes_Rs^{-1}$. As $K$ is torsion free, we deduce from Lemma~\ref{L:torsionfree} that $ys=\sum_{i=1}^{\ell} y_ir_i$. Since $s$ is not a zero divisor
 \[\mathrm{supp} (y)=\mathrm{supp} (ys)\subseteq \mathrm{supp}(\sum_{n\in \N}y_nR).\]
 This finishes the proof of our claim. Now $K\subseteq \bigoplus\limits_{i\in \mathrm{supp} (K)}e_iR$ which is a countably generated direct summand of $F$.
\end{Proof}

\begin{Prop}\label{P:Q-aleph_0noeth} Let $R$ be a ring and let $S\subseteq \Sigma$ satisfy the right Ore condition such that $Q=RS^{-1}$ is right $\aleph_0$-noetherian.
Then the cotorsion pair $(\PP_1, \PP_1^{\perp})$ is of countable
type.
\end{Prop}

\begin{Proof} The result follows by Lemma~\ref{L:countably-generated} using a back and forth argument in
the projective resolution of a module, taking into account that
$\Tor^R_1(M, Q)=0$, for every right $R$-module $M$.
\end{Proof}

\begin{Remark}\label{R:not-countable} \emph {By \cite[Corollary 11]{ST} the cotorsion pair $(\PP_n, \PP_n^\perp)$
is of countable type if and only if every module in $\PP_n$ is $ \PP_n(\mk{\aleph_0})$-filtered.}
\end{Remark}

We show now by an example that the cotorsion pair $(\PP_1,
\PP_1^\perp)$ is not, in general, of countable type and also that
Proposition~\ref{P:Q-aleph_0noeth} cannot be extended to arbitrary
finite projective dimension.

\begin{Ex} \emph{Observe first that if $\m$ is a maximal right ideal
of a ring $R$ then the simple right module $R/\m$ is
$\mk{\aleph_0}$-filtered if and only if $\m \in \mk{\aleph_0}$.}

\emph{1). Let $R$ be the $K$-free algebra generated over the field
$K$ by an uncountable set $X$. Then the two sided ideal generated by
$X$ is an uncountably generated maximal right (or left) ideal $\m$
of $R$. Since, $R$ is a hereditary ring, we infer that the simple
module $R/\m$ has projective dimension 1. In view of
Remark~\ref{R:not-countable}, $(\PP_1, \PP_1^\perp)$ cannot be of
countable type since $R/\m$ is not $P_1(\mk{\aleph_0})$-filtered.}

\emph{2). Let $R$ be a commutative valuation domain such that its
maximal ideal $\m$ is $\aleph _n$-generated. By a result of Osofsky
\cite[Theorem~3.2]{FS}, the projective dimension of $\m$ is $n+1$,
so that the projective dimension of $R/\m$ is $n+2$. If $n>0$ then
$R/\m$ is not $P_{n+2}(\mk{\aleph_0})$-filtered.}
\end{Ex}

\section{Finitistic dimensions of classical rings of quotients}\label{quotients}

We recall the notions of   small and   big finitistic dimension of a
ring $R$. For later  convenience, we introduce also an intermediate
notion.

\begin{Def} The (right) \emph{small finitistic dimension}, f.dim $R$, is the supremum of the projective dimension of the right
$R$-modules in $\PP(\rmod R)$.

The (right) \emph{big finitistic dimension}, F.dim $R$, is the
supremum of the projective dimension of the right $R$-modules in
$\PP$.

We denote by f$_{\aleph_0}$.dim $R$ the supremum of the projective
dimension of the right $R$-modules in $\PP(\mk{\aleph_0})$,

Clearly, f.dim $R\leq $ f$_{\aleph_0}$.dim $R\leq$ F.dim $R$.
\end{Def}

We note the following easy but useful lemma.
\begin{Lemma}\label{L:fin-pres} Let $R$ be a ring and let $C\in \PP_1(\rmod R)$. There is a finitely  generated projective module $P$ and a short exact sequence
\[0\to R^m \to R^n\to C\oplus P\to 0.\]
\end{Lemma}
\begin{Proof} If $C$ is projective, the claim is obvious with $m=0$. Let p.d.$C=1$. By assumption, there exists a short exact sequence $0\to P\to R^k\to C\to 0$ for some $k>0$ and some finitely generated projective module $P$.  Let $P'$ be a projective module such that $P\oplus P'\cong R^m$ for some $m>0$. Then $R^k\oplus P'\oplus P\cong R^{k+m}$ and thus the short exact sequence
\[0\to P\oplus P'\to R^k\oplus P'\oplus P\to C\oplus P\to 0\]
satisfies the  requirements.
\end{Proof}

As before, for a ring $R$, we denote by $\Sigma$ the multiplicative
set of the non zero divisors of $R$.


Let $\C=\{R/rR\mid r \in \Sigma\}\cup R$. Then $\D=\C^{\perp}$ and $\TFcal=\C^{\intercal}$, where $\D$ is the class
of divisible right $R$-modules and $\TFcal$ is the class of torsion
free left $R$-modules.

Clearly $\C\subseteq \PP_1(\rmod R)$.

%

\begin{Prop}\label{P:divisible} Let $R$ be a ring with a classical ring of quotients $Q$. Assume that f.dim $Q$ =0.
Then the following hold.
\begin{enumerate}
\item[(i)] The class $\D$ of divisible right modules is a $1$-tilting class and it  coincides with $\PP_1(\rmod R)^\perp$.
\item[(ii)] The class $\TFcal$ of torsion free left modules 
coincides with $\PP_1(\rmod R)^{\intercal}$.
\end{enumerate}
\end{Prop}

\begin{Proof} Let $C_R\in \PP_1(\rmod R)$. By Lemma ~\ref{L:fin-pres}, $C_R$ fits in a short exact sequence of the form
\begin{equation} \label{eqn:ex-seq}
 0\to R^m \overset  {\mu}\to R^n\to  C\to 0.
\end{equation}
where  the injection $\mu$ can be represented by a $n\times m$ matrix $A$ with entries in
$R$ and acting on the elements of $R^m$ viewed as columns vectors. %
Tensoring the exact sequence (1) by the flat left $R$-module $Q$ we
get the short exact sequence
 \[0\to Q^m\overset {A\otimes 1_Q}\to Q^n\to  C\otimes _RQ\to 0\]
 of right $Q$-modules. Using the assumption f.dim $Q=0$, we conclude that $C\otimes _RQ$ is a projective right $Q$-module. Thus there is a
 splitting map $Q^n\to Q^m$ represented by an $m\times n$ matrix $B'$ with entries in $Q=\Sigma^{-1}R$ such that
 $B'A=I_m$, where $I_m$ is the $m\times m$ identity matrix. Let $r\in \Sigma$ be the product
 of the left denominators of the entries in $B'$, then the matrix $B=rB'$ has entries in $R$, and   $BA=rI_m$.
Thus we have the following commutative diagram:
\[(*)\qquad \xymatrix{R^m\ar[r]^A \ar[dr]_{\overline r }
                               &R^n\ar[d]^B\\
                                 &R^m}\]

 where $\overline r$ denotes the map given by left multiplication by $r$.

(i)  If we show that $\D= \PP_1(\rmod R)^\perp$, then we will have
that $\D$ is a $1$-tilting class, since the cotorsion pair
$(^\perp{}(\PP_1(\rmod R)^\perp), \PP_1(\rmod R)^\perp)$ is a
$1$-tilting cotorsion pair. By definition,  $\D\supseteq \PP_1(\rmod
R)^\perp$. We need to show that $\Ext^1_R(C, D)=0$, for every $C\in
\PP_1(\rmod R)$ and for every $D_R\in \D$.  Applying the functor
$\Hom_R(-, D)$ to  the sequence (1), we obtain the exact sequence
 \begin{equation} 0\to \Hom_R(C, D)\to D^n\stackrel{\small{ \Hom_R(A, D)}}\longrightarrow D^m\to \Ext^1_R(C, D)\to 0
 \end{equation}
 where the map $\Hom_R(A, D)$ is represented by the matrix $A$ acting by right multiplication on elements of $D_R^n$ viewed as row vectors. Applying the functor $\Hom_R(-, D)$ to the commutative diagram (*), we obtain the commutative diagram:

\[\xymatrix{D^m\ar[r]^{\tiny{\Hom(B, D)}} \ar[dr]_{\overline r }
                               &D^n\ar[d]^{\tiny{\Hom(A, D)}}\\
                                 &D^m}\]

 Since the right multiplication by $r$ is surjective on $D$, we conclude that the group homomorphism $\Hom_R(A, D)$ is surjective. Hence, from sequence (2) we infer that $ \Ext^1_R(C, D)=0$.

(ii) By definition,  $\TFcal\supseteq \PP_1(\rmod R)^\intercal$. Let
$Y\in \TFcal$. Applying the functor $-\otimes_R Y$ to sequence (1),
we obtain the exact sequence
 \begin{equation} 0\to \Tor_1^R(C, Y)\to Y^m\overset {A\otimes_R 1_Y}\longrightarrow Y^n\to C\otimes_R Y\to 0 \end{equation}
 where the map ${A\otimes_R 1_Y}$ is represented by the matrix $A$ acting as left multiplication on elements of ${}_RY^m$ viewed as columns vectors. Applying the functor $-\otimes_R Y$ to the commutative diagram (*), we obtain the commutative diagram:

\[\xymatrix{Y^m\ar[r]^{\small{A\otimes_R 1_Y}} \ar[dr]_{\overline r }
                               &Y^n\ar[d]^{\small{B\otimes_R 1_Y}} \\
                                 &Y^m}\]
 Since the left multiplication by $r$ is injective on $Y$, we conclude that the group homomorphism
 ${A\otimes_R 1_Y}$ is injective. Hence, from sequence (3) we infer that $ \Tor^R_1(C,
 Y)=0$. Hence $Y\in \PP_1(\rmod R)^\intercal$ as we wanted to show.
\end{Proof}

In what follows,  $\PP_n(R)$ and $\PP_n(Q)$ will denote the classes
of right modules of projective dimension at most $n$ over $R$ and
$Q$, respectively.

\begin{Lemma}\label{L:tensoring} Let $R$ be a ring with  classical ring of quotients $Q$. Then, a right $Q$-module
$V$ belongs to $ \PP_1(Q)$ if and only if  there is $M_R\in \PP_1(R)$ such that $V=M\otimes_RQ$.
\end{Lemma}
\begin{Proof} The sufficiency is clear. For the only if part, let $V\in \PP_1(Q)$. Without loss of generality we can assume that there is a short exact sequence
\[0\to Q^{(\alpha)}\overset{\mu}\to Q^{(\beta)}\to V\to 0,\]
for some cardinals $\alpha, \beta$.

Let $(d_i\colon i\in \alpha)$ be the canonical basis of the right
$Q$-free module $Q^{(\alpha)}$. The injection $\mu$ is represented
by a column finite matrix $A'$ with entries in $Q=R\Sigma^{-1}$
acting as left multiplication on the basis elements $d_i$. For every
$i\in \alpha$, let $r_i\in \Sigma$ be a common right denominator of
the elements of the $i^{\rm{th}}$-column of $A'$.  Changing the
basis $(d_i\mid i\in \alpha)$ with the basis $(r_id_i\colon i\in
\alpha)$, we can assume that the monomorphism $\mu$ is represented
by a column finite matrix $A$ with entries in $R$. As $R$ is inside
$Q$, we get the short exact sequence
\[0\to R^{(\alpha)}\overset{\nu}\to R^{(\beta)}\to M\to 0,\]
where the map $\nu$ is represented by the matrix $A$. Then it is clear that $M\otimes_RQ\cong V$.
\end{Proof}

A characterization of the rings with  classical ring of quotients $Q$  of big finitistic dimension zero is now
immediate.

\begin{Prop}\label{P:Fin-dim-Q} Let $R$ be a ring with  classical ring of quotients $Q$.  Then, the following are equivalent:
\begin{enumerate}
\item[(i)] For every right $R$-module $M\in \PP_1(R)$, $M\otimes_RQ\in \PP_0(Q)$;
\item[(ii)] F.dim $Q=0$.
\end{enumerate}
\end{Prop}
\begin{Proof}
(i )$\Rightarrow$ (ii).  Assume by way of contradiction that F.dim
$Q>0$. Let $n$ be the least natural number such that there is a non
projective right module $V\in \PP_n(Q)$. Consider a free
presentation $0\to V_1\to Q^{(\alpha)}\to V\to 0$ of $V$. Then
$V_1\in \PP_{n-1}(Q)$, hence $V_1$ is projective. So $V$ has
projective dimension one. By Lemma~\ref{L:tensoring} and condition
(i) we get a contradiction.

(ii) $\Rightarrow$ (i). Obvious because $Q$ is  flat as a left $R$-module.
\end{Proof}

We give now a characterization of rings with classical ring of quotients $Q$ of small finitistic dimension $0$.
\begin{Prop}\label{P:fin-dim-Q} Let $R$ be a ring with  classical ring of quotients $Q$. Then, the following are equivalent:
\begin{enumerate}
\item[(i)] For every right $R$-module $C\in \PP_1(\rmod R)$, $C\otimes_RQ\in \PP_0(\rmod Q)$;
\item[(ii)] f.dim $Q=0$;
\item[(iii)] $\TFcal=(\PP_1(\rmod R))^{\intercal}$;
\item[(iv)]  $\TFcal\supseteq Q \LMod $.
\end{enumerate}
\end{Prop}
\begin{Proof} (i) $\Rightarrow$ (ii). Follows from
Lemma~\ref{L:tensoring} (cf. Proposition~\ref{P:fin-dim-Q}).

(ii )$\Rightarrow$ (iii).  By Proposition~\ref{P:divisible} (ii).

(iii) $\Rightarrow$ (iv). Let $N$ be a left $Q$-module. The left multiplication by an element of $\Sigma$ yields a bijection on $N$. Thus, as a left  $R$-module, $N\in \TFcal$.

(iv )$\Rightarrow$ (i).  Let $C$ be a right $R$-module in $
\PP_1(\rmod R)$ and let $N$ be a left $Q$-module. By hypothesis
$\Tor_1^R(C, N)=0$. As the ring homomorphism $R\to Q$ is an
epimorphism, $0=\Tor_1^R(C, N)\cong\Tor_1^Q(C\otimes_RQ, N)$. So
$C\otimes_RQ$ is a flat right $Q$-module, hence projective, since it
is finitely presented.
\end{Proof}

We consider now a situation which is intermediate between the ones
considered in  Propositions~\ref{P:fin-dim-Q} and
~\ref{P:Fin-dim-Q}.

\begin{Prop}\label{P:omega-dim-Q}  Let $R$ be a ring with  classical ring of quotients $Q$. Then, the following statements
are equivalent:
\begin{enumerate}
\item[(i)] For every right $R$-module $M\in \PP_1(\mk{\aleph_0})$, $M\otimes_RQ\in \PP_0(\qk{\aleph_0})$;
\item[(ii)] f$_{\aleph_0}.$dim $Q=0$;
\item[(iii)] f.dim $Q=0$ and $M\otimes_RQ$ is a pure projective module, for every right $R$-module $M\in \PP_1(\mk{\aleph_0})$;
\item[(iv)] f.dim $Q=0$ and $M\otimes_RQ$ is a Mittag-Leffler right $Q$-module, for every right $R$-module $M\in \PP_1(\mk{\aleph_0})$
\end{enumerate}
\end{Prop}
\begin{Proof}  The equivalence (i) $\Leftrightarrow$ (ii) follows by the definition of  f$_{\aleph_0}.$dim $Q=0$ and by Lemma~\ref{L:tensoring}.

(ii) $\Rightarrow$ (iii). Condition (ii) clearly implies f.dim $Q=0$. Moreover,  for every right $R$-module $M\in \PP_1(\mk{\aleph_0})$, $M\otimes_RQ$ is pure projective right $Q$-module, since by hypothesis it is projective.

(iii )$\Rightarrow$ (i).  Let $M_R\in \PP_1(\mk{\aleph_0})$. Then, as $M_R$ is countably presented and of  projective dimension at most one, it is a direct limit of a countable direct system of the form $(C_n; f_n\colon C_n\to C_{n+1})_{n\in \N}$,  where the right $R$-modules $C_n\in \PP_1(\mk{\aleph_0})$ (\cite[Sec.2]{BH}). Hence $M$ fits in a pure exact sequence of the form
\[0\to \oplus_{n\in \N} C_n\overset{\phi}\to  \oplus_{n\in \N} C_n\to M\to 0\]
where, for every $n \in \N$, $\phi\varepsilon_n=\varepsilon_n
-\varepsilon_{n+1}f_n$ and $\varepsilon_n\colon C_n\to
\oplus_{n\in \N} C_n$ denotes the canonical map. Tensoring by $Q$ we get the pure exact sequence of right $Q$-modules
\[0\to \oplus_{n\in \N} (C_n\otimes_RQ)\overset{\phi\otimes_RQ}\to  \oplus_{n\in \N} (C_n\otimes_RQ)\to M\otimes_RQ\to 0,\]

which is splitting by the hypothesis that $M\otimes_RQ$ is pure projective. Thus $M\otimes_RQ$ is a direct summand of $\oplus_{n\in \N} (C_n\otimes_RQ)$ and for every $n\in \N$, $C_n\otimes_RQ$ is projective right $Q$-module, since f.dim $Q=0$. Thus $M\otimes_RQ$ is projective, too.

(iii) $\Leftrightarrow$ (iv). The equivalence follows by the well known fact that countably generated (hence countably presented) Mittag-Leffler right modules are pure projective \cite[Corollaire 2.2.2 ]{RG}.
\end{Proof}

\section{Orders in rings with finitistic dimension zero}\label{sufficient}

We start by giving a characterization for the equality of the two
classes $\PP_1(\mk{\aleph_0})^\perp$ and $\PP_1(\rmod R)^\perp$.

\begin{Prop}\label{P:suff-cond} Let $R$ be a ring with classical ring of quotients
$Q$ such that f.dim $Q=0$. Then the following statements are
equivalent.

\begin{enumerate}
\item[(i)] f$_{\aleph_0}$.dim $Q=0$
\item[(ii)] Every right $R$-module $M\in \PP_1(\mk{\aleph_0})$ is a summand of a $\PP_1(\rmod R)$-filtered module;
\item[(iii)] the cotorsion pair generated by $\PP_1(\mk{\aleph_0})$ is of finite type.
\end{enumerate}
\end{Prop}
\begin{Proof} Conditions (ii) and (iii)  are equivalent by Fact~\ref{F:Jan-description}.

$(i)\Rightarrow (iii)$. Let $(\A ,\B )$ be the cotorsion pair of
finite type generated by $\PP_1(\rmod R)$. We must show that every
right $R$-module $M$ in $\PP_1(\mk{\aleph_0})$ is in $\A$. As any
module in $\PP _1$ is a direct limit of modules in $\PP_1(\rmod R)$,
by Theorem~\ref{T:Dol-Lid} we only need to show that a right
$R$-module $M$ in $\PP_1(\mk{\aleph_0})$ is Mittag-Leffler with
respect to the class $\PP_1(\rmod R)^{\intercal}$. By
Proposition~\ref{P:fin-dim-Q},  $\PP_1(\rmod R)^{\intercal}$
coincides with the class $\TFcal$ of torsion free left $R$-modules.

We show now that, under our hypothesis, every right $R$-module $M$ in $\PP_1(\mk{\aleph_0})$ is $\TFcal$-Mittag-Leffler.

The assumption  f$_{\aleph_0}$.dim $Q=0$ implies that  $M\otimes_RQ$ is a projective right $Q$-module, hence a
Mittag-Leffler right $Q$-module.

We claim that $M$ is $\Q$-Mittag-Leffler, where $\Q=Q\LMod$.

In fact,  for every right $R$-module $N$ and any left $Q$-module $V$,  $N\otimes_RV\cong
N\otimes_R(Q\otimes_Q V)$. Hence if  $({V_i}; {i\in I})$ is a family of left $Q$-modules, the above remark and the fact that $M\otimes_RQ$ is a projective right $Q$-module imply that the map
\[\rho \colon  M\bigotimes_R \prod _{i\in I}V_i\to \prod _{i\in
I}(M\bigotimes_RV_i)\] is injective.

Let now ${}_RY\in \TFcal$ and consider the exact sequence
\[(1)\quad 0\to R \to Q \to Q/R\to 0.\]
By Lemma~\ref{L:torsionfree}, $\Tor^R_1(
Q/R, Y)=0$. Thus, tensoring by $Y$ the exact sequence (1) we obtain the
embedding
\[(2)\quad 0\to R\otimes_RY \to Q\otimes_RY.\]
Since $Q\otimes_RY$ is a left $Q$-module, Proposition~\ref{P:submodules} implies that $M$ is $\TFcal$-Mittag-Leffler.

$(iii)\Rightarrow (i)$. By Theorem~\ref{T:Dol-Lid}, the cotorsion
pair generated by $\PP_1(\mk{\aleph_0})$ is of finite type then
every module $M\in \PP_1(\mk{\aleph_0})$ is Mittag-Leffler with
respect to the class $(\PP_1(\rmod R))^{\intercal}$. As f.dim $Q=0$,
Proposition~\ref{P:divisible} implies $\TFcal=(\PP_1(\rmod
R))^{\intercal}$.

By Proposition~\ref{P:fin-dim-Q},  $Q \LMod $ is contained in
$\TFcal$. Thus the right $Q$-module $M\otimes_RQ$ is Mittag-Leffler.
The conclusion follows by Proposition~\ref{P:omega-dim-Q}.

\end{Proof}

We now patch together our results in the setting of orders into
$\aleph _0$-noetherian rings.

In the next theorem $\partial$ denotes the Fuchs' divisible module defined in
\cite[VII.1]{FS} for the commutative case and in \cite[\S 5]{hdiv}
for the noncommutative setting. The module  $\partial$ is a
$1$-tilting module generating the cotorsion pair $({}^\perp\D, \D)$
(cf. \cite{F} for the commutative case and
\cite[Proposition~5.5]{hdiv} for the general case).

\begin{Th}\label{P:1} Let $R$ be a ring with an $\aleph_0$-noetherian classical ring of quotients $Q$.
Assume that f.dim $Q=0$. Then   the following statements are
equivalent
\begin{itemize}
\item[(i)] f$_{\aleph_0}$.dim $Q=0$
\item[(ii)] F.dim $Q=0$
\item[(iii)] $(\PP_1, \PP_1^\perp)$ is  of finite type;
\item[(iv)]  Every module of
projective dimension at most one is a direct summand of a
$\PP_1(\rmod R)$-filtered module.
 \item[(v)]  Every module of
projective dimension at most one is a direct summand of a
$\C$-filtered module, where $\C=\{R/rR\mid r\in \Sigma\} \cup
\{R\}$.
\end{itemize}

When the above equivalent statements hold then $(\PP_1,
\PP_1^\perp)=(\PP_1, \D)$ where $\D$ is the class of divisible
modules; so that every divisible module of projective dimension at
most one is a direct summand of a direct sum of copies of
$\partial$. Moreover, every module of projective dimension at most
two is a direct limit of modules in $\PP_2(\rmod R)$.
\end{Th}
\begin{Proof} (i)$\Leftrightarrow$ (ii). Follows from
Fact~\ref{F:countable-type} and Eklof's Lemma
(Fact~\ref{F:Eklof-filt}).

(i)$\Leftrightarrow$ (iii). If  f$_{\aleph_0}$.dim $Q=0$ then, by
Proposition~\ref{P:Q-aleph_0noeth} and
Proposition~\ref{P:suff-cond}, it follows that $(\PP_1,
\PP_1^\perp)$ is of finite type. The converse follows from
Proposition~\ref{P:suff-cond}.

Statements (iii), (iv) and (v) are equivalent by
Fact~\ref{F:Jan-description} and Proposition~\ref{P:divisible}.

When the statements hold then $(\PP_1, \PP_1^\perp)=(\PP_1, \D)$ by
Proposition~\ref{P:divisible}. In this situation, $\partial$ is a
$1$-tilting module generating the cotorsion pair $(\PP_1, \D)$
\cite[Proposition~5.5]{hdiv}. Therefore, by well known results on
tilting cotorsion pairs, $\PP_1\bigcap \D$ is the class $\Add\
\partial$ consisting of direct summands of direct sums of copies of
$\partial$.

The statement on the modules of projective dimension two is a
consequence of Theorem~\ref{T:lim-in P_n}.
\end{Proof}

\section{Orders in semisimple artinian rings and noetherian rings}\label{S:positive}

A semisimple artinian ring has global dimension $0$ and it is
artinian, therefore Theorem~\ref{P:1} applies immediately to orders
into semisimple artinian rings, that is, to semiprime Goldie rings.

\begin{Cor}\label{C:goldie} Let $R$ be a semiprime Goldie ring then the conclusions
of Theorem~\ref{P:1} hold for $R$. In particular, $(\PP_1,
\PP_1^\perp)$ is  of finite type.
\end{Cor}

From the previous Corollary, we single out the case  of commutative
domain, as it completes the results obtained in \cite{L} by
S.~B.~Lee.
 \begin{Cor}\label{C:domains} Let $R$ be a commutative domain then the conclusions
of Theorem~\ref{P:1} hold for $R$. In particular, $(\PP_1,
\PP_1^\perp)$ is  of finite type.
\end{Cor}

Our next goal is to characterize the commutative noetherian rings
such that the cotorsion pair $(\PP _1, \PP _1 ^\perp)$ is of finite
type  as the ones that are orders into artinian rings. Therefore, in
the commutative noetherian case, Theorem~\ref{P:1} gives the best
possible result. We remark however that in
Remark~\ref{R:contravariantly} we will see that the condition f.dim
$Q=0$ is not a necessary condition for the
 cotorsion pair $(\PP_1, \PP_1^\perp)$ to be of finite type.

\begin{Lemma}\label{L:noeth-Q}  Let $R$ be a noetherian commutative ring with classical ring of quotients
$Q$.
Then, f.dim $Q=0$.
\end{Lemma}
\begin{Proof} It is well known that the set of zero divisors of a commutative ring $R$ coincides with
the union of the prime ideals of $R$ associated to $R$. Let $\{P_1, P_2\dots, P_n\}$ be the set of
the prime ideals associated to $R$. For every $1\leq i\leq n$, let $P_iQ$ denote the extension of
$P_i $ in $Q$. Then $\{P_1Q, P_2Q\dots, P_nQ\}$ is the set of prime ideals of $Q$, and by
\cite[Theorem~6.2]{Ms}, it is the set of associated primes of $Q$. Let $P_iQ$ be a maximal ideal
in $Q$ and consider  the localization $Q_{P_iQ}$ of $Q$ at $P_iQ$. Again by
\cite[Theorem~6.2]{Ms}, the maximal ideal of  $Q_{P_iQ}$ is an associated prime of $Q_{P_iQ}$,
hence it consists of zero divisors. This means that the   regular sequences in $Q_{P_iQ}$ are
empty. Hence by the  Auslander Buchsbaum Formula, \cite{AB} or \cite[Theorem~4.4.15]{weibel},
f.dim $Q_{P_iQ}=0$. Since this holds for all maximal ideals of $Q$, we conclude  that any finitely
generated (presented) module of finite projective dimension is flat and, hence, projective.
Therefore, f.dim $Q=0$.
\end{Proof}

\begin{Th}\label{T:noetherian} Let $R$ be a commutative noetherian ring with classical ring of quotients $Q$. Then
the following are equivalent.
\begin{enumerate}
\item[(i)] The cotorsion pair $(\PP_1, \PP_1^\perp)$ is of finite type.
\item[(ii)] F.dim $Q=0$.
\item[(iii)] $Q$ is artinian.
\item[(iv)] the set of  prime ideals associated to $R$ coincides with the set of minimal prime ideals of $R$.
\end{enumerate}
\end{Th}
\begin{Proof} Over any $\aleph_0$-noetherian ring the  cotorsion pair $(\PP_1, \PP_1^\perp)$ is of countable type. Thus for such rings, $(\PP_1, \PP_1^\perp)$ is of finite type if and only if the cotorsion pair generated by $\PP_1(\mk{\aleph_0})$ is of finite type, by Fact~\ref{F:Eklof-filt}.

(i) $\Leftrightarrow$ (ii). By Lemma~\ref{L:noeth-Q}, f.dim $Q=0$.
The above remark and Theorem~\ref{P:1} give the equivalence.

(ii) $\Leftrightarrow$ (iii).  A combination of a result by Bass
\cite{B} and one by Raynaud Gruson \cite{RG} shows that, for a commutative noetherian ring, the
big finitistic dimension equals the Krull dimension. Moreover, a commutative noetherian ring is artinian
if and only if its Krull dimension is zero.

(iii) $\Leftrightarrow$ (iv).  As noted in the proof of
Lemma~\ref{L:noeth-Q}, the prime ideals of $Q$ are exactly the
extension at $Q$ of the associated prime ideals of $R$. Hence the
claim is immediate.
\end{Proof}

Noetherian Cohen-Macaulay rings have an artinian ring of quotients
so they satisfy the above theorem.

\medskip

Kaplansky's characterization of commutative rings with big
finitistic dimension zero (see \cite[pag 1]{B}) combined with
Theorem~\ref{P:1} allows us to prove,

\begin{Remark}\emph{Let $R$ be a commutative ring and $Q$ its total ring
of quotients. Assume that $Q$ is a perfect ring and $\aleph_0$-noetherian.
Then,  $(\PP_1, \PP_1^\perp)$ is of finite type. }
\end{Remark}

\section{Examples}\label{examples}

In this section we exhibit examples and counterexamples for the
finite type of the cotorsion pairs $(\PP_n, \PP_n^\perp)$. Our first
type of examples is based on the following observation.

\begin{Lemma}\label{basic} {\rm Let $R$ be a ring  such that f.dim $R=m<$  F.dim $R$,
then $(\PP_n, \PP_n^{\perp})$ is not of finite type, for all $n>m$.}
\end{Lemma}
\begin{Proof}
By Auslander's Lemma, any direct summand of a $\PP_n(\rmod
R)$-filtered module has projective dimension at most $m$. But, by
assumption, for any $n>m$, $\PP_n(\rmod R)=\PP_m(\rmod R)$ and in
$\PP_n$ there exist modules  of projective dimension greater than
$m$. Therefore, for all $n>m$,  $(\PP_n, \PP_n^{\perp})$ is not of
finite type.
\end{Proof}

In trying to generalize the  results in Section~\ref{S:positive} to
the cotorsion pair $(\PP_n, \PP_n^\perp)$, for $n>1$, the first
thing to keep in mind  are the next two counterexamples showing
that, even over  commutative domains these cotorsion pairs are not
of finite type, in general.

\begin{Exs}\label{E:counterex}{\rm  (i) There is a commutative local noetherian domain such that the cotorsion pair
$(\PP_2, \PP_2^\perp)$ is not of finite type.}

{\rm (ii) If $R$ is a  non Dedekind Pr\"ufer domain, then $(\PP_n,
\PP_n^{\perp})$ is not of finite type, for all $n>1$.}
\end{Exs}
\begin{Proof} An example of the type claimed in (i) is the non Cohen-Macaulay ring
in  \cite[Ex.2.1.18, pag64]{BH}. Let $R=K[[X^4, X^3Y,
XY^3,Y^4]]\subset K[[X,Y]]$, where $K$ is a field and $X, Y$ are
indeterminates. $R$ is a local noetherian domain of  Krull dimension
$2$ and  $X^4, Y^4$ is a system of parameters, but  it is not a
regular sequence. In fact, $Y^4(X^3Y)^2 = X^6Y^6 = X^4(XY^3)^2$, but
$(X^3Y)^2 \notin  (X^4)$, so depth $R=1$. Hence, by
Auslander-Buchsbaum equality  \cite{AB}, f.dim $R=1$ and by
\cite{RG} F.dim $R=2$.  Now the conclusion follows from
Lemma~\ref{basic}.

To prove (ii) recall that  finitely presented modules over a
Pr\"ufer domain $R$ have projective dimension at most one, hence
$\PP_n(\rmod R)=\PP_1(\rmod R)$, for every $n\geq 1$. Now our
statement will follow from Lemma~\ref{basic}, once we have proved
that  in a non Dedekind Pr\"ufer domain $\PP _1\subsetneq \PP _2$.

To this aim note that a   non Dedekind Pr\"ufer domain is a non
noetherian ring, hence it has a countably generated ideal $I$ that
is not finitely generated. Being $R$ semihereditary, $I$ is flat,
and, since $R$ is a domain, it is countably presented. As  $I$ is
flat and countably presented it has projective dimension at most
$1$. Since in  a domain the projective ideals are finitely
generated,  we deduce that $I$ has projective dimension exactly $1$.
Therefore $R/I\in \PP _2\setminus \PP _1$.
\end{Proof}

On the positive side, we consider the case of an Iwanaga-Gorenstein
ring, that is a  left and right noetherian ring $R$ such that the
right module $R_R$ has finite injective   dimension and the left
module ${}_RR$ has also finite injective  dimension. In this case,
both dimensions coincide. The ring $R$ is said to be an
$n$-Iwanaga-Gorenstein if these dimensions are both $n$.

 \begin{Ex}\label{E:n-Iwanaga}  {\rm If $R$ is an  $n$-Iwanaga-Gorenstein ring, then $(\PP_n, \PP_n^\perp)$
 is of finite type.}
\end{Ex}
\begin{Proof}
 It was shown in \cite[Theorem~3.2]{AHT} that if $R$ is an $n$-Iwanaga- Gorenstein ring,
 then f.dim $R=$ F.dim $R=n$ and that the cotorsion pair generated by $\PP(\rmod R)$ is the
 $n$-tilting cotorsion pair corresponding to the $n$-tilting module
 $T=\bigoplus\limits_{0\leq i\leq n} I_i$ where $0\to R\to I_0\to I_1\to \dots \to I_n\to 0$ is
 a minimal
 injective coresolution of $R$.  Moreover, in \cite{AHT} it is shown that $\PP=  ^\perp{}(T^\perp)$.
 Hence, $(\PP, \PP^\perp)= (\PP_n, \PP_n^\perp)$ is of finite type.
 \end{Proof}

 \begin{Ex}{\rm If $R$ is a commutative   Gorenstein ring then it is Cohen-Macaulay.
Hence, by Theorem~\ref{T:noetherian} the cotorsion pair   $(\PP_1,
\PP_1^\perp)$ is always of finite type and it is generated by
$\{R/rR\mid \mbox{ $r$ regular element of $R$}\}$.

 If $R$ is $n$-Gorenstein, we do not know whether $(\PP_m, \PP_m^\perp)$ is of
finite type for  $1<m<n$, cf. Proposition~\ref{P:Birge}.}
\end{Ex}

 \begin{Ex}{\rm
   (i)  If f.dim $R=0$, then $\PP_n(\rmod R)= \PP_0(\rmod R)$, for every $n$. Hence, $(\PP_n, \PP_n^\perp)$ is of finite type if and only if F.dim $R=0.$

(ii) If $R$ is a right noetherian ring, right self-injective, then all right projective modules are injective. Hence F.dim $R=0$ and so  for every $n\in \N$, $(\PP_n, \PP_n^\perp)= (\PP_0, \PP_0^\perp)$ is of finite type.
}
 \end{Ex}

\bigskip

Next we consider the case of an artin algebra, that is a finitely generated algebra over a commutative artin ring.

Recall that  a
subclass $\X$ of $\PP(\rmod R)$ is said to be \emph{contravariantly
finite} if every $M\in \rmod R$  admits an
$\X$-precover (cover), that is there exist $X\in \X$ and a morphism $f\colon X\to M$ such that $\Hom_R(X', X)\to \Hom_R(X', M)$ is surjective for every $X'\in \X$.

Auslander and Reiten \cite{AR} proved  a fundamental result, namely
that  if $\PP(\rmod R)$ is contravariantly finite, then the small
finitistic dimension of $R$ is finite.

 Huisgen-Zimmermann and Smal\o\ in \cite{HS} strengthened  Auslander-Reiten's result by
 proving that, if $\PP(\rmod R)$ is contravariantly finite, then the big finitistic dimension of
 $R$ coincides with its small finitistic dimension.

 In \cite[Theorem 4.3]{AT} Angeleri and Trlifaj showed that, for any
right noetherian ring  $R$,  f.dim $R\leq n$ if and only if the cotorsion pair generated by
$\PP(\rmod R)$ is an $n$-tilting cotorsion pair. Moreover, they prove that for an artin algebra
$R$,  $\PP(\rmod R)$ is contravariantly finite in $\rmod R$ if and only if the tilting module
corresponding to the cotorsion pair generated by $\PP(\rmod R)$ can be taken to be finitely
generated.
 Thus, as a consequence of all these results we have:

 \begin{Ex}\label{E:AT}{\rm \label{E:controvariantly finite} Let $R$ be an artin algebra.
 Assume that $\PP(\rmod R)$ is contravariantly finite in $\rmod R$. Let
 f.dim $R= n$(= F.dim $R$). Then, $(\PP_n, \PP_n^\perp)$ is of finite type.}
\end{Ex}

\begin{Proof} By the preceding remarks and  \cite[Corollary 3.6]{AT}.
 \end{Proof}

\begin{Remark}\label{R:contravariantly}
\emph{ In contrast with our previous discussion on rings with
classical ring of quotients with finitistic dimension $0$, we note
that an artin algebra coincides with its classical ring of
quotients. So Example~\ref{E:controvariantly finite} shows that
there exists a ring  with classical ring of quotients of small
finitistic dimension greater than zero such that the cotorsion pair
$(\PP_1, \PP_1^\perp)$ is of finite type.}
\end{Remark}

Since over right perfect rings, direct limits of module of finite
projective dimension $n$ are still of finite projective dimension
$n$, we have the following general observation.

\begin{Prop}\label{P:right-perfect} Let $R$ be a right perfect  ring.
Assume that f.dim $R=n$ and F.dim $R>n$, for some $n\geq 1$. Then the cotorsion pair $(\PP_n, \PP_n^\perp)$ is
not of finite type.
\end{Prop}
\begin{Proof} By hypothesis, there exists a right module $M$ of projective dimension exactly $n+1$.
Assume, by way of contradiction that  $(\PP_n, \PP_n^\perp)$ is of
finite type. By Theorem~\ref{T:lim-in P_n}, $M_R$ is a direct limit
of objects in $\PP_{n+1}(\rmod R)$ which coincides with $\PP_n(\rmod
R)$, by assumption. Since $R$ is right perfect, p.d.$M\leq n$ (see
\cite[Theorem P]{B}), a contradiction.
\end{Proof}

In \cite{S}
 Smal\o\  constructs a family of examples of finite dimensional algebras $R_n$, such that
 f.dim $R_n=1$ an F.dim $R_n=n$ for every $n\in \N$. So that, for $n>1$,
 $R_n$ satisfies the hypothesis of
 Proposition~\ref{P:right-perfect}.

\begin{Ex} {\rm   In \cite{IST}, Igusa, Smal\o\ and Todorov construct an example of a finite dimensional monomial
algebra such that f.dim $R=1=$ F.dim $R$.  However, as proved in
\cite[Sec. 5]{AT2},  $(\PP_1, \PP_1^{\perp})$ is not of finite
type.}
\end{Ex}

We devote the rest of the section to give an example showing that
the finite type property of $(\PP_n, \PP_n^\perp)$ is not
inherited, in general, by $(\PP_{n-1}, \PP_{n-1}^\perp)$. We recall
that this was mentioned in the second statement of
Remark~\ref{R:asc-desc}.

As the example will be a quotient of a path algebra, we find it
more convenient to think the modules as representations of the
associated quiver. So from now on our statements will involve left modules.

We will examine
the behavior of the functor $\mathrm{Ext}$ with respect to inverse
limits of modules over artin algebras.

To this aim recall that if we have a (countable) inverse system
$(H_n)_{n\in \N}$ and a sequence of morphisms
\[\dots\to H_{n+1}\stackrel{h _n}{\to} H_n\to\dots\to H_3\stackrel{h_2}{\to}
H_2\stackrel{h_1}\to H_1\] then $\varprojlim H_n$ fits into the
exact sequence
\[0\to \varprojlim H_n\to \prod _{n\in \N}H_n\stackrel{\Delta}\to
\prod _{n\in \N}H_n\] where $\Delta =(\mathrm{Id_{H_n}}-h_n)_{n\in
\N }$. By definition, $\mathrm{coker} \Delta =\varprojlim ^1
(H_n)_{n\in \N}$. Moreover, if the inverse system
$(H_n)_{n\in \N}$ satisfies the Mittag-Leffler condition, then $\varprojlim ^1
(H_n)_{n\in \N}=0$. (See \cite[\S3.5]{weibel}.

A result similar to the next one appears in  \cite[\S 3]{CB} with a different approach.

\begin{Lemma}\label{L:alginverse-limit} Let $R$ an artin algebra. Let $(M_n, f_n\colon M_{n+1}\to
M_{n})_{n\in \N}$ be an inverse system of finitely generated left
$R$-modules. Then, for any left $R$-module $A$ and for any $k\ge 0
$, $\mathrm{Ext}_R^k (A,\varprojlim M_n)\cong \varprojlim
\mathrm{Ext}_R^k (A,M_n)$.
\end{Lemma}

\begin{Proof} The ring $R$ has a duality that we denote by $D$,
and any finitely generated module $M$ satisfies that $M\cong
D(D(M))$. Also $\varprojlim M_n\cong D(\varinjlim D(M_n))$. This
allows us to conclude that, being dual modules, $M_n$ and
$\varprojlim M_n$ are pure injective. As for any pure injective
module $E$, any direct system $(A_\alpha, u_{\alpha \,
\beta}\colon A_\alpha \to A_\beta)_{\alpha \le \beta \in I}$ and any $k\ge 0$ there is an
isomorphism
\[\mathrm{Ext}_R^k (\varinjlim A_\alpha,E)\cong \varprojlim
\mathrm{Ext}_R^k (A_\alpha ,E)\] to prove our claim we may assume
that $A$ is finitely generated. Moreover, since all the syzygies of  a finitely generated
module are again finitely generated, by dimension shifting, it is enough to show the result for the case $k=1$.

We shall use repeatedly that a countable inverse system of
finitely generated modules over an artin ring satisfies the Mittag-Leffler
condition.

Set $M=\varprojlim M_n$. Using the canonical presentation of the inverse limit, we have an exact sequence:
\[0\to M\to \prod _{n\in \N} M_n\overset{\Delta}{\to}\prod _{n\in \N}
M_n \to 0,\]
Applying the functor $\Hom_R(A,- )$ to it we obtain the canonical presentation of the
inverse limit of the Mittag-Leffler inverse system $(\Hom_R (A,  M_n), \Hom_R (A,
f_n))_{n\in \N}$, hence we get the exact sequence:
\[0\to \Hom_R (A,M )\to \prod _{n\in \N} \Hom_R (A,  M_n)\overset{\Delta _H}{\to}
\prod _{n\in \N} \Hom_R (A,  M_n)  \to 0 \]
Therefore the following sequence is also exact
\[0\to \Ext_R^1(A,M)\to \prod _{n\in \N}\Ext_R
^1(A,M_n)\overset{\Delta _E}{\to}\prod _{n\in \N}\Ext_R
^1(A,M_n).\] As $\Delta _E$ is the canonical map of the
presentation of the inverse limit of the inverse system $(\Hom_R
(A,  M_n), \Hom_R (A, f_n))_{n\in \N}$ we deduce that
$\Ext_R^1(A,M)\cong \varprojlim \Ext_R^1(A,M_n)$.

\end{Proof}

\begin{Cor}\label{C:inverse-limitp1} Let $R$ be an  artin algebra. Let $(M_n, f_n\colon
M_{n+1}\to M_{n})_{n\in \N}$ be an inverse system of finitely
generated left $R$-modules. If, for any $n\in
\N$, $M_n\in \PP_1(R \lmod)^\perp$ then $\varprojlim M_n\in \PP_1(R
\lMod)^\perp$.
\end{Cor}
\begin{Proof} By Lemma~\ref{L:alginverse-limit}, $\varprojlim M_n\in\PP_1(R
\lmod)^\perp$. The conclusion follows by the same argument as in the first part of the proof of Lemma~\ref{L:alginverse-limit}, since any module of projective dimension at most 1 is a direct limit of finitely presented modules of projective dimension at most 1 and the module $\varprojlim M_n$ is pure injective.
\end{Proof}

\begin{Ex}\label{E:Birge}  {\rm [Communicated by B. Huisgen-Zimmermann]\\
Consider the quiver $Q$ given by}

$$\xymatrix{ 3 \ar[r]^\delta & 1 \ar@/^/[r]^\alpha \ar@/_/[r]^\beta &
2 \ar@/^2pc/[l]_\gamma\ar[r]^\epsilon & 4
 }$$

{\rm  Let $K$ be a field and consider the path algebra $R = KQ/I$
where the ideal $I$ is generated by: $\epsilon\beta$, $\gamma\beta$,
$\beta\delta$, $\epsilon\alpha\delta$; all paths leaving the vertex
1 that have length at least 3; all paths leaving the vertex 2 that
have length at least 2. Then, the following hold:}
 \begin{enumerate}

 \item{\rm  By \cite{H} and \cite{HS2}, $\PP(R \lmod )$ is contravariantly finite and f.dim $R=$F. dim $R=2$, so every module in $\PP_2$ is a direct limit of objects in
$\PP_2(R \lmod)$.}

 \item{\rm  By \cite{HS2}, $\PP_1(R \lmod)$ fails to be contravariantly finite.}
\end{enumerate}
\end{Ex}

\begin{Prop}\label{P:Birge} Let $R$ be the finite dimensional algebra defined in Example~\ref{E:Birge}. Then
$(\PP_2(R$-$\mathrm{Mod}), \PP_2(R$-$\mathrm{Mod})^\perp)$ is of
finite type, but $(\PP_1(R$-$\mathrm{Mod}),
\PP_1(R$-$\mathrm{Mod})^\perp)$ fails to be of finite type.
\end{Prop}

\begin{Proof}
For $i\in \{1,2,3,4\}$, let $P_i=R e_i$ denote the indecomposable
projective left modules of $R$ and let $I_i=E(S_i)$ denote the
indecomposable injective left modules.

Let $J=P_3\oplus R \epsilon \alpha \oplus R \gamma \alpha \oplus R
\epsilon \oplus R \gamma$. Note that $J$ is a two-sided ideal of $R$
and that $R /J$ is isomorphic to the \emph{Kronecker algebra} that
we shall denote by $\Lambda$. The left $\Lambda$ modules are left $R
$ modules via the projection $R \to R /J=\Lambda $.

Consider the simple regular modules over $\Lambda $:
\[ V_\lambda =\xymatrix{  K \ar@/^/[r]^\lambda \ar@/_/[r]^1 &
K
 } \quad\mbox{for every $\lambda \in K$}; \qquad   V_\infty = \xymatrix{  K \ar@/^/[r]^1 \ar@/_/[r]^0 &
K.
 }\]

Then,\\
(i) For every $\lambda \in K$, $V_\lambda$ is a finitely generated
$R$-module of projective dimension $1$.\\
(ii) $V_\infty \in \PP_1(R \lmod )^\perp$.

In fact,  as an $R$-module, $V_\lambda \cong P_1/R
(\alpha -\lambda \beta)$ and $R (\alpha -\lambda \beta)\cong
P_2$. Therefore  (i)  holds.
To verify (ii), note that $V_\infty$ is a quotient of $I_4$ and recall that $ \PP_1(R \lmod )^\perp$
contains the injective modules and is closed under epimorphic images.

For any $\lambda \in K$, denote by $T_\lambda$ the corresponding
$\Lambda $-Pr\"ufer module and by $t_\lambda$ the corresponding
tube in $\Lambda \lmod$. As $T_\lambda$ and the modules
in $t_\lambda$ are filtered by $V_\lambda$, condition (i) above tells us
 that they are modules in $\A ={}^\perp(
\PP_1(R \lmod)^\perp)$.

As $\B= \PP_1(R \lmod)^\perp$ is a tilting class, it is closed by
direct limits and extensions. Hence, by condition (ii) above, all
the modules in $t_\infty$ and the Pr\"ufer module $T_\infty$ are
in $\B$. By Corollary~\ref{C:inverse-limitp1}, we can also
conclude that the adic module $Z_\infty$ is in $\B$. Therefore,
for any set $I$, $Z_\infty^{(I)}\in \B$.

Now we are ready to proceed as in \cite{AT2} to conclude that
$(\PP_1(R$-$\mathrm{Mod}), \PP_1(R$-$\mathrm{Mod})^\perp)$ is not of
finite type.

By \cite[Proposition~3]{R}, if $T_\lambda$ is any  of the Pr\"ufer
modules of the Kronecker algebra, then the generic module $Q$ is a
direct summand of $T_\lambda^{\N}$. Since for finite dimensional
algebras, $\PP_1$ is closed under products, taking $\lambda \in K$
we deduce that the generic module $Q$ has projective dimension $1$
viewed as  an $R$-module.  Since $Z_{\infty}$ is the dual of a
Pr\"ufer module it is pure injective, however it is not $\Sigma$-pure
injective.  By results due to Okoh \cite[Proposition 1 and
Remark]{O},  $\Ext_{\Lambda}
^1(Q,Z_{\infty}^{(\N)})= 0$ would imply $Z_{\infty}^{(\N)}$ pure injective. We conclude that
$\Ext_R^1(Q,Z_{\infty}^{(\N)})\neq 0$ and therefore, by
Proposition~\ref{P:finite}, $(\PP_1, \PP_1^\perp)$ is not of
finite type, since $\PP_1^\perp\neq \B$.
\end{Proof}

\end{document}